\documentclass[aos,seceqn,nameyear,dvips]{arximspdf}

\doi{10.1214/10-AOS803}
\volume{38}
\issue{5}
\pubyear{2010}
\firstpage{2751}
\lastpage{2780}

\makeatletter

\newtheorem{theorem}{Theorem}
\newtheorem{lemma}{Lemma}
\newproclaim{definition}{Definition}

\makeatother

\begin{document}
\begin{frontmatter}

\title{On universal oracle inequalities related to high-dimensional
linear models}
\runtitle{On universal oracle inequalities}

\begin{aug}
\author[A]{\fnms{Yuri} \snm{Golubev}\corref{}\thanksref{t1}\ead[label=e1]{golubev.yuri@gmail.com}}
\runauthor{Y. Golubev}
\affiliation{CNRS, Universit\'e de Provence}
\address[A]{
CMI\\
Universit\'e de Provence\\
39 rue F. Joliot-Curie\\
13453 Marseille\\
France\\
\printead{e1}} 
\end{aug}

\thankstext{t1}{Supported in part by ANR Grant 07-0234-01.}

\received{\smonth{8} \syear{2009}}
\revised{\smonth{1} \syear{2010}}

%
\begin{abstract}
This paper deals with recovering an unknown vector $\theta$ from the
noisy data $Y=A\theta+\sigma\xi$, where $A$ is a known $(m\times n)$-matrix
and $\xi$ is a white Gaussian noise. It is assumed that $n$ is
large and $A$ may be severely ill-posed. Therefore, in order to
estimate $\theta$, a spectral regularization method is used, and our
goal is to choose its regularization parameter with the help of the
data $Y$. For spectral regularization methods related to the so-called
ordered smoothers [see \citeauthor{K} \textit{Ann. Statist.} \textbf{22} (\citeyear{K}) 835--866], we propose new penalties in the
principle of empirical risk minimization. The heuristical idea
behind these penalties is related to balancing excess risks. Based
on this approach, we derive a sharp oracle inequality controlling
the mean square risks of data-driven spectral regularization methods.
\end{abstract}

%

\end{frontmatter}

\section{Introduction and main results}\label{s1}
In this paper, we consider a classical problem of
recovering an unknown vector $\theta=(\theta(1),\ldots,\theta
(n))^{\top}\in\mathbb{R}^n$ in the standard linear model
%
%
\begin{equation}\label{equ.0}
Y=A\theta+\sigma\xi,
\end{equation}
where $A$ is a known $(m\times n)$-matrix and $\xi=(\xi(1),\ldots
,\xi(m))^\top
$ is a standard white Gaussian noise in $\mathbb{R}^m$ with
$\mathbf{E}\xi(k)=0, \mathbf{E}\xi^2(k)=1, k=1,\ldots,m$. The
noise level $\sigma$ in (\ref{equ.0}) is assumed to be known.

We start out by considering the maximum likelihood estimate of $\theta$
\[
\hat{\theta}_0={\mathop{\arg\min}_{\theta\in\mathbb{R}^n}}\|
Y-A\theta\|^2,
\]
where $ \|x\|^2=\sum_{k=1}^mx^2(k).
$ It is easily seen that $
\hat{\theta}_0= (A^{\top}A)^{-1}A^{\top}Y $ and that the mean
square risk
of this estimator is computed as follows:
%
%
\begin{eqnarray}\label{equ.1}
\mathbf{E} \|\hat\theta_0-\theta\|^2 &=&\sigma^2\mathbf{E}
\|(A^{\top}A)^{-1}A^{\top}\xi\|^2
=\sigma^2 \operatorname{trace}[(A^{\top}A)^{-1}] \nonumber\\[-8pt]\\[-8pt]
&=&\sigma^2\sum_{k=1}^n
\lambda^{-1}(k),\nonumber
\end{eqnarray}
where $\lambda(k)$ and $\psi_k\in\mathbb{R}^n$ are eigenvalues
and eigenvectors of $A^{\top}A$
\[
A^{\top}A\psi_k =\lambda(k)\psi_k.
\]

In this paper, it is assumed solely that
$\lambda(1) \ge\lambda(2) \ge\cdots\ge\lambda(n)$. So, $A$ may
be severely ill-posed and
(\ref{equ.1}) reveals the principal difficulty in $\hat{\theta}_0$:
\textit{its risk may be very large when $n$ is large or when $A$ has
a large condition number.}

The simplest way to improve $\hat\theta_0$ is to suppress
large $\lambda^{-1}(k)$ in (\ref{equ.1}) with the help of a linear
smoother; that is, to estimate $\theta$ by $H\hat\theta_0$, where
$H$ is a properly chosen $(n\times n)$-matrix. In what follows,
we deal with smoothing matrices admitting the following representation
$H=H_\alpha(A^{\top}A)$,
where $H_\alpha(\lambda)$ is a function $\mathbb{R}^+\rightarrow
[0,1]$ which depends on a regularization parameter $\alpha\in[0,\bar
\alpha]$ such that
\[
\lim_{\alpha\rightarrow0}H_\alpha(\lambda)=1,\qquad \lim_{\lambda
\rightarrow
0}H_\alpha(\lambda)=0.
\]

This method is called \textit{spectral
regularization} [see \citet{EHN}] since $A^\top A$ and $H_\alpha
(A^{\top}A)$ have the same eigenvectors.
Summarizing, we estimate $\theta$ with the help of the following
family of linear estimators
\[
\hat\theta_\alpha=
H_\alpha(A^{\top}A)(A^{\top}A)^{-1}A^{\top}Y
\]
and our main goal is to choose the best estimator within this\vspace*{1pt} family,
or equivalently, the best regularization parameter $\alpha$. Note that
$\alpha$ controls the mean square risk
of~$\hat\theta_\alpha$,
%
%
\begin{equation} \label{equ.2}
L_\alpha(\theta) \stackrel{\mathrm{def}}{=} \mathbf{E}\|\hat
\theta_\alpha-\theta\|^2 =\sum_{k=1}^n[1-
h_\alpha(k)]^2\langle\theta,\psi_k\rangle^2+\sigma^2
\sum_{k=1}^n \lambda^{-1}(k) h_\alpha^2(k),\hspace*{-25pt}
\end{equation}
where here and below we denote for brevity
\[
h_\alpha(k) \stackrel{\mathrm{def}}{=} H_\alpha[\lambda(k)]
\quad\mbox{and}\quad
\langle\theta,\psi_k\rangle\stackrel{\mathrm{def}}{=} \sum_{l=1}^n
\theta(l)\psi_k(l).
\]

According to (\ref{equ.2}), the variance of $\hat\theta_\alpha$ is
always smaller than that of the maximum likelihood estimate, but $\hat
\theta_\alpha$ has a nonzero bias and adjusting properly $\alpha$ we
may improve $\hat\theta_0$.
Note that this improvement may be significant
if $\langle\theta,\psi_k\rangle^2$ are small for
large $k$.

In practice, a good choice of $H_\alpha(\cdot)$ is a delicate
problem related to the numerical complexity of $\hat\theta_\alpha$. For
instance, to make use of the spectral cut-off with
$H_\alpha(\lambda) =\mathbf{1}\{ \lambda\ge\alpha\}$, one has to compute
the singular value decomposition (SVD) of $A$. For very large $n$, this
numerical problem may be difficult or even infeasible.

The very popular Tikhonov--Phillips [see, e.g., \citet{TA}] regularization
\[
\hat\theta_{\alpha}=\mathop{\arg\min}_{\theta
} \{\|Y-A\theta\|^2+\alpha\|\theta\|^2 \}
\]
does not require SVD. In this case, $\hat\theta_\alpha$ is computed
as a root of the linear equation
\[
(\alpha I+ A^{\top}A) \hat\theta_\alpha= A^{\top}Y
\]
and therefore
$ H_\alpha(\lambda)={\lambda}/(\lambda+\alpha)$. It is worth
pointing out
that this regularization technique is good solely for ill-posed $A$.

Another widespread regularization technique is due to \citet{L}.
This method is based on a very simple idea: to compute recursively a
root of equation
\[
A^{\top}A \theta= A^{\top}Y.
\]
Since $
A^{\top}Y=[A^{\top}A-aI] \theta+a \theta$ for all $a>0$,
we get $ \theta=[I-a^{-1}A^{\top}A]\theta
+a^{-1}A^{\top}Y. $ This formula motivates Landweber's iterations
defined by
\[
\hat{\theta}_{k}=[I-a^{-1}A^{\top}A]\hat{\theta}_{k-1}
+a^{-1}A^{\top}Y.
\]
Thus, we can estimate $\theta$ without computing SVD and without
solving linear equations.
It is easily seen that these iterations converge if $\lambda(1)<a$ and
that the corresponding spectral regularization function is given by
%
%
\begin{equation}\label{equ.3}
H_k(\lambda)=1- \biggl(1-\frac{\lambda}{a} \biggr)^{k+1}.
\end{equation}
The regularization parameter of the Landweber method is usually defined
by $\alpha=1/k$.
Note that in spite of its iterative character, the numerical complexity of
the Landweber method may be hight. Indeed, when the noise is very
small, $H_k(\lambda)$ should be close to $1$, and (\ref{equ.3})
implies that
\[
k\gtrsim\operatorname{cond}(A)\stackrel{\mathrm{def}}{=}\frac{\lambda
(1)}{\lambda(n)}.
\]
This means that if $A$ is severely ill-posed, the number of iterations
may be
very large, thus making the method infeasible. A substantial
improvement of Landweber's iterations is provided by the so-called
$\nu$-method [see, e.g., \citet{EHN} and \citet{BHMR}].

All the above-mentioned regularization methods are particular cases of
the so-called ordered smoothers
[see \citet{K}] defined as follows.
\begin{definition}\label{d1}
The family of sequences $\{h_\alpha(k), \alpha\in(0,\bar{\alpha}],
k \in\mathbb{N}^+\}$ is
called ordered smoother if:
\begin{enumerate}
\item For any given $\alpha\in(0,\bar{\alpha}]$, $ h_\alpha(k)\dvtx
\mathbb{N}^+\rightarrow[0,1]$ is a monotone function of $k$.
\item If for some $ \alpha_1,\alpha_2\in(0,\bar{\alpha}]$ and $
k'\in\mathbb{N}^+$,
$h_{\alpha_1}(k') < h_{\alpha_2}(k')$, then $
h_{\alpha_1}(k) \le h_{\alpha_2}(k)$ for all
$k\in\mathbb{N}^+$.
\end{enumerate}
\end{definition}

It was Kneip who noted that from a probabilistic viewpoint, all ordered
smoothers are equivalent to
the spectral cut-off with $h_\alpha(k)=\mathbf{1}\{\lambda(k) \ge
\alpha\}$.
This profound fact plays an essential role in adaptive estimation since
it helps to analyze precisely statistical risks of feasible data-driven
regularization methods. This is why in this paper we deal solely with
the ordered smoothers.

Whatever inversion method is used, the principal question usually
arising in practice is how
to choose its regularization parameter. Traditional theoretic approach
to this problem is related to the minimax theory; see, for example,
\citet{MR} and \citet{S}. However, this approach provides the
smoothing
parameters depending strongly on an a priory information about $\theta
$ which is hardly available in practice. The only one way to improve
this drawback is to use data-driven regularizations. In statistical
literature, one can find several general approaches for constructing
such methods. We cite here, for instance,
the Lepski method which has been adopted to inverse problems in
\citet
{Ma}, \citet{BH}, \citet{BHMR} and the model selection
technique which
was implemented in \citet{LL}.

In this paper, we take the classical way related to the famous
principle of unbiased risk estimation which goes back to \citet{A}.
The heuristical motivation of this approach is based on the idea that a
good data-driven regularization should minimize in some sense the risk
$L_\alpha(\theta)$ [see
(\ref{equ.2})].
This idea is put into practice with
the help of the
empirical risk minimization suggesting to compute data-driven
regularization parameters as follows:
%
%
\begin{equation}\label{equ.4}
\hat\alpha=\mathop{\arg\min}_{\alpha\in(0,\bar\alpha]}
{R}_{\alpha}[Y,\operatorname{Pen}],
\end{equation}
where
\[
{R}_\alpha[Y,\operatorname{Pen}]=
\|\hat\theta_0-\hat\theta_\alpha\|^2+\sigma^2\operatorname{Pen}(\alpha),
\]
and $\operatorname{Pen}(\alpha)\dvtx(0,\bar{\alpha}]\rightarrow\mathbb{R}^+$ is a
given penalty function.
The most important problem in this approach is related to the choice of
the penalty.
Intuitively, we want that the method mimics the oracle smoothing parameter
$\alpha^*=\arg\min_{\alpha}L_\alpha(\theta)$. This is why we are
looking for a minimal penalty that ensures the following inequality:
%
%
\begin{equation}\label{equ.5}
L_\alpha(\theta)\lesssim
{R}_{\alpha}[Y,\operatorname{Pen}]+\mathcal{C},
\end{equation}
where $\mathcal{C}$ is a random variable that does not depend on
$\alpha$.
It is easily seen that in the considered statistical model,
\[
\mathcal{C}=-\|\theta-\hat\theta_0\|^2=-\sigma^2\sum_{k=1}^n
\lambda^{-1}(k)\xi^2(k).
\]

Traditional approach to solving (\ref{equ.5}) is based on the unbiased
risk estimation defining the penalty as a root of the equation
\[
L_\alpha(\theta)=
\mathbf{E}{R}_{\alpha}[Y,\operatorname{Pen}]+\mathbf{E}\mathcal{C}.
\]
Unfortunately, in spite of its very natural motivation, this penalty is
not good for ill-posed problems [see \citet{CavG} for more details].

The main idea in this paper is to compute the penalty in a little bit
different way, namely as a minimal function assuring the following
inequality:
%
%
\begin{equation}\label{equ.6}\qquad
\mathbf{E}\sup_{ \alpha\le\bar{\alpha}} \bigl[L_\alpha(\theta)-
{R}_{\alpha}[Y,\operatorname{Pen}]-\mathcal{C} \bigr]_+ \le K\mathbf{E} \bigl[ L_{\bar
{\alpha}}(\theta)-
{R}_{\bar{\alpha}}[Y,\operatorname{Pen}]-\mathcal{C} \bigr]_+,
\end{equation}
where $[x]_+=\max\{0,x\}$ and
$K>1$ is a constant. The heuristical motivation behind this approach is
rather transparent: we are looking for a minimal penalty that balances
all excess risks uniformly in $\alpha\in(0,\bar{\alpha}]$. Recall
that the excess risk is defined as the difference between the risk of
the estimate and its empirical risk. Note that according to (\ref
{equ.5}), we may focus on the positive part of the excess risk, and
that equation (\ref{equ.6}) guarantees that for any data driven
smoothing parameter $\hat\alpha$
\[
\mathbf{E} \bigl[L_{\hat\alpha}(\theta)-
{R}_{\hat\alpha}[Y,\operatorname{Pen}]-\mathcal{C} \bigr]_+\le K \mathbf{E} \bigl[ L_{\bar
{\alpha}}(\theta)-
{R}_{\bar{\alpha}}[Y,\operatorname{Pen}]-\mathcal{C} \bigr]_+.
\]

In order to explain how one can compute good penalties assuring (\ref
{equ.6}), we begin with the spectral representation of the underlying
statistical problem.
We can check easily that
\[
y(k)\stackrel{\mathrm{def}}{=}\langle Y,\psi_k\rangle\lambda^{-1}(k)=
\langle\theta,\psi_k\rangle+\sigma\lambda^{-1/2}(k) \xi(k),
\]
where $\xi(k)$ are i.i.d. $\mathcal{N}(0,1)$. With this notation,
$\hat\theta_\alpha$ admits the following representation:
\[
\langle\hat\theta_\alpha,\psi_k\rangle=h_\alpha(k)y(k)=h_\alpha
(k)\theta(k)+\sigma h_\alpha(k)\lambda^{-1/2}(k)\xi(k),
\]
where $\theta(k)=\langle\theta,\psi_k\rangle$, and
%
%
\begin{eqnarray} \label{equ.7}
\|\hat\theta_0-\hat\theta_\alpha\|^2&=&\sum_{k=1}^n[1-h_\alpha(k)]^2
y^2(k),\nonumber\\[-8pt]\\[-8pt]
\|\theta-\hat\theta_\alpha\|^2&=&\sum_{k=1}^n[\theta(k)-h_\alpha
(k)y(k)]^2.\nonumber
\end{eqnarray}

In what follows, it is assumed that the penalty has the following
structure:
\[
\operatorname{Pen}(\alpha) =2 \sum_{k=1}^n
\lambda^{-1}(k) h_\alpha(k)+(1+\gamma)Q(\alpha),
\]
where $\gamma$ is a positive number and $Q(\alpha), \alpha>0$, is a
positive function to be defined later on.
Then the excess risk is computed as follows:
%
%
\begin{eqnarray}\label{equ.8}\quad
&& L_\alpha(\theta)-{R}_{\alpha}[Y,\operatorname{Pen}]-\mathcal{C}\nonumber\\
&&\qquad =\sigma^2 \sum_{k=1}^n \lambda^{-1}(k)[2h_\alpha(k)-h_\alpha
^2(k)]\bigl(\xi^2(k)-1\bigr)-(1+\gamma)
\sigma^2Q(\alpha)\\
&&\qquad\quad{} -2\sigma\sum_{k=1}^n \lambda^{-1/2}(k)[1-h_\alpha(k)]^2\xi
(k)\theta(k).\nonumber
\end{eqnarray}

Our first idea in solving (\ref{equ.6}) is to use the fact that the
absolute value of the cross term
\[
2\sigma\sum_{k=1}^n \lambda^{-1/2}(k)[1-h_\alpha(k)]^2\xi
(k)\theta(k)
\]
is typically smaller than $L_\alpha(\theta)$ (for more details, see
Lemma \ref{lemma.9} below). Therefore, omitting this term in (\ref
{equ.6}), we get the following inequality for
$Q(\alpha)$:
%
%
\begin{equation}\label{equ.9}
\mathbf{E}\sup_{\alpha\le\bar{\alpha}}[\eta_\alpha-(1+\gamma
)Q(\alpha)]_+ \le K\mathbf{E}[\eta_{\bar{\alpha}}-(1+\gamma
)Q(\bar{\alpha})]_+ ,
\end{equation}
where
\[
\eta_\alpha=\sum_{k=1}^n \lambda^{-1}(k)[2h_\alpha(k)-h_\alpha
^2(k)]\bigl(\xi^2(k)-1\bigr).
\]
Usually, computing the minimal function $Q(\alpha)$ assuring (\ref
{equ.9}) is a hard numerical problem. However, when
$h_\alpha(k)$ is a family of ordered smoothers it can be solved
relatively easily.
The main idea is to find a feasible solution $Q^\circ(\alpha)$ of the
marginal inequality
%
%
\begin{equation}\label{equ.10}
\mathbf{E}[\eta_\alpha-Q^\circ(\alpha)]_+\le\mathbf{E}[\eta
_{\bar{\alpha}}-Q^\circ(\bar\alpha)]_+
\end{equation}
and then to show that $(1+\gamma)Q^\circ(\alpha)$ satisfies (\ref{equ.9}).
To solve (\ref{equ.10}), we use the following inequality:
%
%
\begin{equation}\label{equ.11}
\mathbf{E}[\eta-x]_+^p\le\Gamma(p+1)\lambda^{-p}\exp(-\lambda x)
\mathbf{E}\exp(\lambda\eta),
\end{equation}
which holds for any random variable $\eta$ and for any $\lambda>0$.
Its proof follows from the Chernoff bound.
Without loss of generality, we may assume that $Q^\circ(\bar\alpha
)=0$. Therefore, according to the Cauchy--Schwarz inequality
\[
\mathbf{E}[\eta_{\bar{\alpha}}-Q^\circ(\bar\alpha)]_+ \le\sqrt
{\mathbf{E}\eta_{\bar{\alpha}}^2}= D(\bar\alpha),
\]
where
\[
D(\alpha)=\sqrt{\mathbf{E}\eta^2_\alpha}= \Biggl\{2\sum_{k=1}^n \lambda
^{-2}(k)[2h_\alpha(k)-h_\alpha^2(k)]^2 \Biggr\}^{1/2}.
\]
%
Hence, $Q^\circ(\alpha)$ is computed as a root of equation
\[
\inf_{\lambda}\exp[-\lambda Q^\circ(\alpha)]\mathbf{E}\exp
(\lambda\eta_\alpha)=D(\bar{\alpha}).
\]

It is not difficult to check with a little algebra that
%
%
\begin{equation}\label{equ.12}
Q^\circ(\alpha)=2D(\alpha)\mu_\alpha\sum_{k=1}^n\frac{\rho
^2_\alpha(k)}{1-2\mu_\alpha\rho_\alpha(k)},
\end{equation}
where $\mu_\alpha$ is a root of equation
%
%
\begin{equation}\label{equ.13}
\sum_{k=1}^n F[\mu_\alpha\rho_\alpha(k)]=\log\frac{D(\alpha
)}{D(\bar{\alpha})}
\end{equation}
and
%
%
\begin{eqnarray}\label{equ.14}
F(x)&=&\frac{1}{2}\log(1-2x)+x +\frac{2x^2}{1-2x}\nonumber\\[-8pt]\\[-8pt]
\rho_\alpha(k)&=&\sqrt{2}D^{-1}(\alpha)\lambda^{-1}(k)[2h_\alpha(k)
-h_\alpha^2(k)].\nonumber
\end{eqnarray}

The only one numerical difficulty in computing $Q^\circ(\alpha)$ is
related to
(\ref{equ.13}). However note that in the proof of Lemma \ref{lemma.7}
it is shown that
\[
f(\mu)=\sum_{k=1}^n F[\mu\rho_\alpha(k)]
\]
is a strictly monotone function and therefore (\ref{equ.13}) may be
solved exponentially fast. Note also that Lemma \ref{lemma.7} provides
lower and upper bounds for $Q^{\circ}(\alpha)$ and $\mu_\alpha$. In
particular, for some constant $C>0$
%
%
\begin{equation}\label{equ.67}
C^{-1}D(\alpha)\log^{1/2} \frac{D(\alpha)}{D(\bar\alpha)}\le
Q^{\circ}(\alpha)\le C D(\alpha)\log\frac{D(\alpha)}{D(\bar
\alpha)}.
\end{equation}

The next theorem shows that $Q^\circ(\alpha)$ computed as a root of
the marginal inequality (\ref{equ.10}) satisfies the global inequality
(\ref{equ.9}).
\begin{theorem} \label{theorem.1} Let $Q^\circ(\alpha)$ be defined
by (\ref{equ.12})--(\ref{equ.14}). Then for any $\gamma>r\ge0$,
\[
\mathbf{E}\sup_{\alpha\le\bar{\alpha}} [\eta_\alpha-(1+\gamma
)Q^{\circ}(\alpha) ]_+^{1+r}\le\frac{CD^{1+r}(\bar{\alpha
})}{(\gamma-r)^3},
\]
where here and throughout the paper $C$ denotes a generic constant.
\end{theorem}

The following theorem represents the main result in this paper. It
controls the performance of the empirical risk minimization by the
penalized oracle risk defined by
\[
r(\theta)\stackrel{\mathrm{def}}{=}\inf_{\alpha\le\bar{\alpha}} \bar
{R}_{\alpha}(\theta),
\]
where
\[
\bar{R}_{\alpha}(\theta)\stackrel{\mathrm{def}}{=}\mathbf{E}\{
{R}_{\alpha}[Y,\operatorname{Pen}]+\mathcal{C}\} =L_\alpha(\theta) +(1+\gamma
)\sigma^2Q^{\circ}(\alpha).
\]
\begin{theorem}\label{theorem.2} Let
$
\operatorname{Pen}(\alpha) =2 \sum_{k=1}^n
\lambda^{-1}(k) h_\alpha(k)+(1+\gamma)Q^\circ(\alpha)
$ with $Q^\circ(\alpha)$ defined by (\ref{equ.12})--(\ref{equ.14}).
Then the mean square risk of $\hat\theta_{\hat\alpha}$ with the
data-driven smoothing parameter $\hat\alpha$ defined by (\ref
{equ.4}) satisfies the following upper bound:
%
%
\begin{equation}\label{equ.15}\quad
\mathbf{E}\|\theta-\hat\theta_{\hat\alpha}\|^2\le
r(\theta) \biggl\{1 +
\biggl[\frac{C}{\gamma} \log^{-1/2}\frac{Cr(\theta)}{\sigma^2 D(\bar
{\alpha})}+\frac{ C\sigma^2D(\bar{\alpha})}{\gamma^4r(\theta)} \biggr]^{1/2}
\biggr\},
\end{equation}
which holds true uniformly in $\theta\in\mathbb{R}^n$.
\end{theorem}

Below, we discuss briefly some statistical aspects of this theorem.
\begin{enumerate}
\item Equation (\ref{equ.15}) represents a particular form of the
so-called oracle inequality
\[
\mathbf{E}\|\theta-\hat\theta_{\hat\alpha}\|^2\le
r(\theta) \biggl\{1+\Psi\biggl[\frac{\sigma^2 D(\bar{\alpha})}{r(\theta)} \biggr]
\biggr\},
\]
where $\Psi(\cdot)$ is a bounded function such that $\lim
_{x\rightarrow0}\Psi(x)=0$. This means that if the ratio $\sigma^2
D(\bar{\alpha})/r(\theta)$ is small, then the risk of the method is
close to the risk of the penalized oracle. On the other hand, if this
ratio is large, then the risk of the method is of order of the oracle risk.

Note also that (\ref{equ.15}) is a universal oracle inequality which
holds true whatever is the ill-posedness of the underlying inverse
problem. It generalizes the corresponding oracle inequalities in
\citet
{CavG} and \citet{G} obtained for the spectral cut-off method.

\item Theorem \ref{theorem.2} reveals some difficulties related to the
data-driven choice of the regularization parameter in the
Tikhonov--Phillips method. Recalling that for this method $H_\alpha
(\lambda)=\lambda/(\alpha+\lambda)$, we obtain
\begin{eqnarray*}
D^2(\alpha)&=&2\sum_{k=1}^n \lambda^{-2}(k)[2h_\alpha(k)-h_\alpha
^2(k)]^2\ge
2\sum_{k=1}^n \lambda^{-2}(k) h_\alpha^2(k)\\&=&
2\sum_{k=1}^n \frac{1}{[\alpha+\lambda(k)]^2}\ge\frac{2n}{[\alpha
+\lambda(n)]^2}.
\end{eqnarray*}
Since $Q^\circ(\alpha)\ge D(\alpha)\asymp\sqrt{n}$, it is clear
that the penalized oracle risk of the Tikhonov--Phillips regularization
may be very large compared to the risk of the method computed for given
$\alpha$. This means that in practice, the Tikhonov--Phillips
regularization with a data-driven smoothing parameter may fail.

Note however, that this drawback can be easily improved with the help
of high-order Tikhonov--Phillips regularizations computed
as follows:
\[
\hat\theta_\alpha^{(k+1)}=\mathop{\arg\min}_\theta\bigl\{\bigl\|A\hat
\theta^{(k)}_\alpha-A\theta\bigr\|^2+\alpha\|\theta\|^2 \bigr\},
\]
where $\hat\theta^{(1)}_\alpha$ stands for the standard
Tikhonov--Phillips regularization.
One can check with a little algebra that the corresponding smoothers
are given by
$H_\alpha^{(k)}(\lambda)=\lambda^{k}(\alpha+\lambda)^{-k}, k\ge
2$, and everything goes smoothly in this case.

\item If the inverse problem is not severely ill-posed; that is,
$\lambda(k)\ge Ck^{-\beta}$ for some $\beta>0$, then, according to
(\ref{equ.67}), for reasonable spectral regularizations
%
%
\begin{equation}\label{equ.a}
\sum_{k=1}^n h_\alpha^2(k)\lambda^{-1}(k) \gg Q^{\circ}(\alpha),
\end{equation}
when $\alpha$ is small.
This means that the risk of the penalized oracle is close to the risk
of the ideal oracle $\inf_{\alpha\le\bar{\alpha}}L_\alpha(\theta)$.

This remark together with the famous \citet{P} minimax theorem shows
that our method results in adaptive asymptotically (as $\sigma
\rightarrow0$) minimax regularizations. To demonstrate this, suppose
for simplicity that $n=\infty$ and that $\theta$ belongs to the
following ellipsoidal body
\[
\Theta(W) = \Biggl\{\theta\in l_2(1,\infty)\dvtx\sum_{k=1}^\infty\langle
\theta,\psi_k\rangle^2 b^2[\lambda^{-1}(k)]\le W \Biggr\},
\]
where $|b(x)|$ is a nondecreasing function such that for some $p,q\in
(0,\infty)$,
$
Cx^{p} \le b(x)\le Cx^q.
$
Then it follows from \citet{P} that as $\sigma\rightarrow0$
\begin{eqnarray*}
{\inf_{\hat\theta}\sup_{\theta\in\Theta(W)}}\|\theta-\hat\theta
\|^2&=&\bigl(1+o(1)\bigr)\inf_{h}\sup_{\theta\in\Theta(W)}L[\theta,h]\\
&=&\bigl(1+o(1)\bigr)\sup_{\theta\in\Theta(W)}\inf_{h}L[\theta
,h]\\
&=&\bigl(1+o(1)\bigr)\sup_{\theta\in\Theta(W)}\inf_{\alpha\in(0,\bar
{\alpha}]}L[\theta,h^*_\alpha],
\end{eqnarray*}
where $\inf$ at the left-hand side is taken over all estimators,
\[
L[\theta,h]=
\sum_{k=1}^n[1-
h(k)]^2\langle\theta,\psi_k\rangle^2+\sigma^2
\sum_{k=1}^n \lambda^{-1}(k) h^2(k)
\]
and
$
h_\alpha^*(k)=[1-\alpha| b[\lambda^{-1}(k)]|]_+.
$
Recall that from a statistical viewpoint, the main drawback in this
minimax result is that the optimal smoothing parameter
\[
\alpha^*(W) =\mathop{\arg\min}_{\alpha\in(0,\bar{\alpha}]} \sup
_{\theta\in\Theta(W)}L[\theta,h^*_\alpha]
\]
depends on the size $W$ of $\Theta(W)$ which is hardly known in
practice. In order to overcome this difficulty, one may use the
data-driven regularization $\hat\alpha$ with $h_\alpha(k)=h^*_\alpha(k)$.
Noticing that this family of smoothers is ordered, we get according to
Theorem \ref{theorem.2} and (\ref{equ.a})
\[
{\sup_{\theta\in\Theta(W)}}\|\theta-\hat\theta_{\hat\alpha}\|
^2=\bigl(1+o(1)\bigr)\inf_{\hat\theta}\sup_{\theta\in\Theta(W)}\|\theta
-\hat\theta\|^2\qquad \mbox{as } \sigma\rightarrow0.
\]

Another interesting situation is related to the case when the inverse
problem is severely ill-posed; that is, when the eigenvalues of
$AA^\top$ are exponentially decreasing, $\lambda(k)\approx\exp
(-\beta k)$ with some $\beta>0$. Then for small $\alpha$
\[
\sum_{k=1}^n h_\alpha^2(k)\lambda^{-1}(k) \ll Q^{\circ}(\alpha).
\]
This means that the risk of the penalized oracle is essentially greater
than that of the ideal oracle. In this situation, Theorem \ref
{theorem.2} provides an upper bound similar to \citet{G}. It is worth
pointing out that neither (\ref{equ.15}) nor the extra penalty can be
improved in this case [for more details, see \citet{G}].
\end{enumerate}

\section{Proofs}\label{sec2}

\subsection{An exponential chaining inequality}\label{sec21}
Let $\xi_t$ be a separable zero mean random process on $ \mathbb{R}^+$.
Denote for brevity
\[
\Delta_\xi(t_1,t_2)=\xi_{t_1}-\xi_{t_2}.
\]
We begin with a general fact similar to Dudley's entropy bound [see,
e.g., \citet{VW}].
\begin{lemma}\label{lemma.1} Let $\sigma^2_u, u\in\mathbb{R}^+$, be
a continuous strictly increasing function with $\sigma^2_0=0$. Then
for any $\lambda>0$,
%
%
\begin{eqnarray}\label{equ.16}
&&\log\mathbf{E}\exp\biggl\{\lambda\max_{0< u\le
t}\frac{\Delta_\xi(u,t)}{\sigma_t} \biggr\}\nonumber\\[-8pt]\\[-8pt]
&&\qquad\le\frac{\log(2)\sqrt
{2}}{\sqrt{2}-1}
+ \max_{0< u< v\le t}\max_{| z|\le\sqrt{2}/(\sqrt{2}-1)}\log
\mathbf{E}\exp\biggl\{z\lambda
\frac{\Delta_\xi(u,v)}{\bar{\Delta}_\sigma(v,u)} \biggr\},\nonumber
\end{eqnarray}
where
$
\bar{\Delta}_\sigma(v,u)=\sqrt{|\sigma^2_v-\sigma^2_u|}.
$
\end{lemma}
\begin{pf}
The proof of (\ref{equ.16}) is based on the standard chaining argument
[for more details, see \citet{VW}]. Denote for brevity by $t_-(B)$ and
$t_+(B)$ left and right elements of a closed subset $B$ in $\mathbb{R}^+$.
First, we construct a dyadic partition of $[0,t]$.
Let
\[
\mathcal{T}_1^1 = \biggl\{u\ge0\dvtx\sigma^2_u\le\frac{\sigma^2_t}{2}
\biggr\},\qquad
\mathcal{T}_2^1 = \biggl\{u\le t\dvtx\sigma^2_u> \frac{\sigma^2_t}{2} \biggr\}.
\]
Next, we partition $\mathcal{T}_1^1$ and $\mathcal{T}_2^1$ as follows:
\begin{eqnarray*}
\mathcal{T}_1^2 &=& \biggl\{u\in\mathcal{T}_1^1\dvtx\sigma^2_u\le\frac
{\sigma^2_{t_+(\mathcal{T}_1^1)}}{2} \biggr\},\\
\mathcal{T}_2^2 &=& \biggl\{u\in
\mathcal{T}_1^1\dvtx\sigma^2_{u}> \frac{\sigma^2_{t_+(\mathcal
{T}_1^1)}}{2} \biggr\},\\
\mathcal{T}_3^2 &=& \biggl\{u\in\mathcal{T}_2^1\dvtx\sigma^2_{u}\le\frac{
\bar{\Delta}^2_\sigma[t_+(\mathcal{T}_2^1),t_-(\mathcal
{T}_2^1)]}{2} \biggr\},\\
\mathcal{T}_4^2 &=& \biggl\{u\in\mathcal{T}_2^1\dvtx\sigma
^2_{u}> \frac{ \bar{\Delta}^2_\sigma[t_+(\mathcal
{T}_2^1),t_-(\mathcal{T}_2^1)]}{2} \biggr\}.
\end{eqnarray*}
Doing so, after $p$ steps, we get partitions $\mathcal{T}_j^{p}, j=1,
\ldots, 2^p$, such that for any $x,y\in\mathcal{T}_j^k$
%
%
\begin{equation}\label{equ.17}
\bar{\Delta}^2_\sigma(y,x)\le2^{-k}\sigma^2_t.
\end{equation}
With the sets $\mathcal{T}^k_j, j=1,\ldots,2^k$ we associate the set
of their right points
\[
\tau^k=\bigcup_{j=1}^{2^k}\{ t_+(\mathcal{T}^k_j)\},
\]
and for any point $x\in\tau^k $ we denote by $\tau_{k-1}(x)$ the
nearest point
in $\tau^{k-1}$. So, by~(\ref{equ.17}), for any $v\in\tau^{p}$
%
%
\begin{equation}\label{equ.18}
\bar{\Delta}^2_\sigma(\tau_{p-1}(v),v)\le2^{-p+1}\sigma^2_{t}.
\end{equation}

With this notation, for any $u\in\tau^p$, setting $\tau_0(v)=t$,
we obtain
%
%
\begin{eqnarray}\label{equ.19}
\xi_u-\xi_t&=&\sum_{k=1}^{p} \bigl[\xi_{\tau_{k}(u)}-\xi_{\tau_{k-1}(u)}\bigr]
\le\sum_{k=1}^{p} \sup_{v\in\tau^{k}}\bigl[\xi_v-\xi_{\tau
_{k-1}(v)}\bigr]\nonumber\\[-8pt]\\[-8pt]
&=&\sum_{k=1}^{p} \sup_{v\in\tau^{k}}\bar{\Delta}_\sigma[v,\tau_{k-1}(v)]
\times\frac{\Delta_\xi[v,\tau_{k-1}(v)]}{\bar{\Delta}_\sigma
[v,\tau_{k-1}(v)]}.\nonumber
\end{eqnarray}

To bound the right-hand side at the above display, we use the
elementary inequality
%
%
\begin{equation}\label{equ.20}
\log\mathbf{E} \exp\biggl[ \sum_{k} q(k) \eta(k) \biggr] \le
\sum_{k} q(k) \log\mathbf{E} \exp[\eta(k)],
\end{equation}
which holds for any random variables $ \eta(k) $ and any given $
q(k)\ge0 $ with
$\sum_{k} q(k) = 1 $. The proof of (\ref{equ.20}) follows immediately
from the convexity of $
\exp(x) $ which implies
\begin{eqnarray*}
&&\mathbf{E} \exp\biggl\{
\sum_{k} q(k)
\{ \eta(k)- \log\mathbf{E} \exp[\eta(k)] \}
\biggr\} \\
&&\qquad \le\biggl\{
\sum_{k} q(k)
\mathbf{E} \exp\bigl[ \eta(k)- \log\mathbf{E} \exp[\eta(k)] \bigr]
\biggr\} = 1.
\end{eqnarray*}

Applying (\ref{equ.20}) with
\begin{eqnarray*}
q(k)&=&\frac{2^{-k/2}}{ \sum_{l=1}^{p} 2^{-l/2}},\\
\eta(k)&=&\frac{\lambda}{\sigma_t q(k)}\sup_{v\in\tau^{k}}\bar
{\Delta}_\sigma[v,\tau_{k-1}(v)]
\times\frac{\Delta_\xi[v,\tau_{k-1}(v)]}{\bar{\Delta}_\sigma
[v,\tau_{k-1}(v)]},
\end{eqnarray*}
we obtain by (\ref{equ.18}) and (\ref{equ.19})
%
%
\begin{equation}\label{equ.21}
\log\mathbf{E}\exp\biggl\{\lambda\sup_{ u\in\tau^{p}
}\frac{\Delta_\xi(u,t)}{\sigma_t} \biggr\} \le\sum_{k=1}^p \lambda(k)
\log\mathbf{E}
\exp[\eta(k)].
\end{equation}

It is easily seen that for any $\lambda>0$
%
%
\begin{eqnarray}\label{equ.22}
&&\mathbf{E}\exp[\eta(k)]\nonumber\\
&&\qquad =\mathbf{E}\exp\biggl\{\frac{\lambda}{q(k)\sigma_t}
\sup_{v\in\tau^{k}}\bar{\Delta}_\sigma[\tau_{k-1}(v),v]
\times\frac{\Delta_\xi[v,\tau_{k-1}(v)]}{\bar{\Delta}_\sigma
[\tau_{k-1}(v),v]} \biggr\}\nonumber\\[-8pt]\\[-8pt]
&&\qquad \le\sum_{v\in\tau^{k}}\mathbf{E}
\exp\biggl\{ \frac{\lambda}{q(k) \sigma_t}\bar{\Delta}_\sigma[\tau_{k-1}(v),v]
\times\frac{\Delta_\xi[v,\tau_{k-1}(v)]}{\bar{\Delta}_\sigma
[\tau_{k-1}(v),v]}
\biggr\}\nonumber\\
&&\qquad \le2^{k}\sup_{u< v}\sup_{|z|\le\sqrt{2}/(\sqrt{2}-1)}\mathbf{E}
\exp\biggl\{\lambda z \frac{\Delta_\xi(u,v)}{\bar{\Delta}_\sigma(v,u)}
\biggr\}.\nonumber
\end{eqnarray}
In the above equation, it was used that
$
\sum_{s=1}^\infty2^{-s/2}=1/(\sqrt{2}-1)
$.

Finally, substituting (\ref{equ.22}) into (\ref{equ.21}), we arrive at
\begin{eqnarray*}
&&\log\mathbf{E}\exp\biggl\{\lambda\sup_{u\in\tau^{p}
}\frac{\Delta_\xi(u,t)}{\sigma_{t}} \biggr\}
\\
&&\qquad\le \log(2)\sum_{k=1}^p
k q(k) +
\sup_{u< v}\sup_{| z|\le\sqrt{2}/(\sqrt{2}-1)}\mathbf{E}
\exp\biggl\{\lambda z \frac{\Delta_\xi(u,v)}{\bar{\Delta}_\sigma(v,u)}
\biggr\},
\end{eqnarray*}
thus completing the proof.
\end{pf}

\subsection{Ordered processes}\label{sec22}
\begin{definition}
A zero mean process $\xi_t,t\in\mathbb{R}^+$ is called ordered if
there exists a continuous strictly monotone scaling function
$\sigma^2_t, t\in\mathbb{R}^+$ and some $\Lambda> 0$ such that
%
%
\begin{equation}\label{equ.23}
\sup_{u,v\in\mathbb{R}^+\dvtx u \neq v}\mathbf{E}\exp\biggl[ \Lambda\frac
{\Delta_\xi(u,v)}{\bar{\Delta}_\sigma(v,u)} \biggr]<\infty.
\end{equation}
\end{definition}

A banal example of an ordered process is a standard Wiener process $w_t$.
In this case, $\sigma^2_t=t$ and obviously
\[
\mathbf{E}\exp\biggl[ \lambda\frac{\Delta_w(u,v)}{\sqrt{|v-u|}} \biggr]=\exp
(\lambda^2/2).
\]
\begin{lemma}\label{lemma.2}
Let $\xi_t$ be an ordered process with $\xi_0=0$. Then there exists a constant
$C$ such that for all $1\le p,q\le2$, uniformly in $z>0$,
%
%
\begin{equation} \label{equ.24}
\mathbf{E}\sup_{t\ge0}[\xi_t-z\sigma^q_t]_+^p \le
\frac{C}{z^{p/(q-1)}},
\end{equation}
where $[x]_+=\max(0,x)$.
\end{lemma}
\begin{pf}
Without loss of generality, we may assume that
that
$\lim_{t\rightarrow\infty}\sigma^2_t=\infty$. For any integer
$k\ge0$, define $t_k(z)$ as a root of the equation
\[
\sigma^{q-1}_{ t_k(z)}=\frac{2^{1/(q-1)} pk}{\Lambda z}.
\]
%
Then we
have
%
%
\begin{eqnarray}\label{equ.25}\qquad
\mathbf{E}\sup_{t\ge0}[\xi_t-z\sigma^q_t]_+^p &\le& \sum
_{k=0}^\infty\mathbf{E}\sup_{t\in[t_k(\gamma),t_{k+1}(z)]}
[\xi_t-z\sigma^q_t]_+^p
\nonumber\\
&\le&
\sum_{k=0}^\infty\mathbf{E}\sup_{t\in[t_k(z),t_{k+1}(z)]}
[\xi_t-z\sigma^q_{t_k}]_+^p\\
&\le&
{\mathbf{E}\sup_{0\le t \le
t_{1}(z)}}|\xi_t|^p+ \sum_{k=1}^\infty
\mathbf{E} \Bigl[\sup_{0\le t\le
t_{k+1}(z)}\xi_t- z\sigma^q_{t_k(z)} \Bigr]_+^p.\nonumber
\end{eqnarray}
According to Lemma \ref{lemma.1}, the first term at the right-hand
side of the
above inequality is bounded as follows:
%
%
\begin{equation}\label{equ.26}
{\mathbf{E}\sup_{0\le t \le t_{1}(z)}}|\xi_t|^p\le Cp^p
\sigma^p_{t_1(z)}=\frac{Cp^p} {z^{p/(q-1)}},
\end{equation}
whereas the second one, by (\ref{equ.11}) and Lemma \ref{lemma.1}, is
controlled by
\begin{eqnarray*}
&&\mathbf{E} \Bigl[\sup_{ t\le
t_{k+1}(\gamma)}\xi_t- z\sigma^q_{t_k(z)} \Bigr]_+^p
\\
&&\qquad=
\sum_{k=1}^\infty\sigma^{p}_{t_{k+1}(z)} \mathbf{E}
\biggl[\sup_{ t\le t_{k+1}(z)}\frac{\xi_t}{\sigma_{t_{k+1}(z)}}\ge\frac
{z\sigma^q_{t_k(z)}}{\sigma_{t_{k+1}(z)}} \biggr]^p_+\\
&&\qquad\le C\sum_{k=1}^\infty
\sigma^{p}_{t_{k+1}(z)}\exp\biggl[-\Lambda\frac{z\sigma^q_{t_k(z)}}
{\sigma_{t_{k+1}(z)}} \biggr]\\
&&\qquad\le \frac{C}{z^{p/(q-1)}\lambda
^{1/(q-1)}}\sum_{k=1}^\infty{(k+1)^{p/(q-1)}} \exp\biggl[-\frac
{2^{1/(q-1) }pk}{(1+1/k)^{1/(q-1)}} \biggr]\\
&&\qquad\le\frac{C}{z^{p/(q-1)}}.
\end{eqnarray*}
So, combining the above inequality with
(\ref{equ.25}) and (\ref{equ.26}), we arrive at (\ref{equ.24}).
\end{pf}

The next very simple lemma is useful for understanding the fact that
the ordered process is controlled by its variance $\sigma_t^2$.
\begin{lemma}\label{lemma.3} Let $\xi_t,t\in\mathbb{R}^+$, be a
random process
such that
%
%
\begin{equation} \label{equ.27}
\mathbf{E}\sup_{t\ge0}[\xi_t-z\sigma^q_t]_+^p \le
\frac{C}{z^{p/(q-1)}},
\end{equation}
for any $z>0$ and some $p\ge1$, $q>1$.
Then there exists a constant $C'$ such that for any random variable
$\tau\in\mathbb{R}^+$
\[
[\mathbf{E}|\xi_\tau|^p]^{1/p}\le C'[\mathbf{E}\sigma^{qp}_\tau]^{1/(pq)}.
\]
\end{lemma}
\begin{pf}
According to (\ref{equ.27}) and Minkowski's inequality, we obviously have
\begin{eqnarray*}
[\mathbf{E}|\xi_\tau|^p]^{1/p}&\le&
\{\mathbf{E}|\xi_\tau-z\sigma_\tau^q+z\sigma_\tau^q|^p\}^{1/p}\\
&\le&
\Bigl\{\mathbf{E}\max_t[\xi_t-z\sigma_t^q]^p\Bigr\}^{1/p}
+z[\mathbf{E}\sigma_\tau^{pq}]^{1/p}\\
&\le&
z[\mathbf{E}\sigma_\tau^{pq}]^{1/p}+\frac{C}{z^{1/(q-1)}}
\end{eqnarray*}
and minimizing the right-hand side in $z$ we finish the proof.
\end{pf}

\subsection{Ordered processes related to the spectral regularization}\label{sec23}
In this section, we focus on typical ordered processes related to the
empirical risk minimization.

For given $\alpha^\circ\in(0,\bar\alpha]$,
define the following Gaussian processes:
\begin{eqnarray*}
\xi^{+}_\alpha&=&\sum_{k=1}^n
[h_{\alpha^\circ}(k)-h_{\alpha^\circ+\alpha}(k)]b(k) \xi(k),\qquad 0<
\alpha\le\bar\alpha-\alpha^\circ, \\
\xi^{-}_\alpha&=&\sum_{k=1}^n
[h_{\alpha^\circ}(k)-h_{\alpha^\circ-\alpha}(k)]b(k) \xi(k),\qquad 0\le
\alpha\le
\alpha^\circ,
\end{eqnarray*}
where $\xi(k)$ are i.i.d. $\mathcal{N}(0,1)$ and $ \sum_{k=1}^n
b^2(k)< \infty. $ It is easily seen that $\xi^+_{\alpha}$ and $\xi
^-_{\alpha}$
are
ordered processes. Indeed, since they are Gaussian, we can choose
\[
\sigma_{\alpha}^{\pm} =\sqrt{\mathbf{E}(\xi^{\pm}_{\alpha})^{2}}
\]
and it suffices to check that
\[
|\mathbf{E}(\xi^{\pm}_{\alpha_1})^{2}-\mathbf{E}(\xi^{\pm
}_{\alpha_2})^{2}|\ge
\mathbf{E}(\xi^{\pm}_{\alpha_1}-\xi^{\pm}_{\alpha_2})^{2},
\]
or equivalently,
\[
\mathbf{E}\xi^{\pm}_{\alpha_1} \xi^{\pm}_{\alpha_2}\ge\min\{
(\sigma^{\pm}_{\alpha_1})^2,(\sigma^{\pm}_{\alpha_2})^2\}.
\]
If $\alpha_1\le\alpha_2$, then we have
\begin{eqnarray*}
\mathbf{E}[\xi^{\pm}_{\alpha_1}]^2&=& \sum_{k=1}^n
[h_{\alpha^\circ}(k)-h_{\alpha^\circ\pm\alpha_1}(k)]
[h_{\alpha^\circ}(k)-h_{\alpha^\circ\pm\alpha_1}(k)]b^2(k)\\
&\le&\sum_{k=1}^n
[h_{\alpha^\circ}(k)-h_{\alpha^\circ\pm\alpha_1}(k)]
[h_{\alpha^\circ}(k)-h_{\alpha^\circ\pm\alpha_2}(k)]b^2(k)
\\
&=&\mathbf{E}\xi^\pm_{\alpha_1}\xi^\pm_{\alpha_2}.
\end{eqnarray*}

Therefore, according to Lemma \ref{lemma.2}, we get
\begin{eqnarray*}
\mathbf{E}\sup_{0\le\alpha\le
\bar\alpha-\alpha^\circ} \Biggl[\xi^+_{\alpha}-z\Biggl(\sum_{k=1}^n
[h_{\alpha^\circ}(k)-h_{\alpha^\circ+\alpha}(k)]^2 b^2(k) \Biggr)^{q/2} \Biggr]_+^p
&\le&\frac{C(p,q)}{z^{p/(q-1)}}, \\
\mathbf{E}\sup_{0\le\alpha\le
\alpha^\circ} \Biggl[\xi^-_{\alpha}-z\Biggl(\sum_{k=1}^n
[h_{\alpha^\circ}(k)-h_{\alpha^\circ-\alpha}(k)]^2 b^2(k) \Biggr)^{q/2} \Biggr]_+^p
&\le&\frac{C(p,q)}{z^{p/(q-1)}},
\end{eqnarray*}
and combining these inequalities we arrive at the following lemma.
\begin{lemma}\label{lemma.4}
For any $z>0$,
\begin{eqnarray*}
&&\mathbf{E}\sup_{0\le\alpha\le\bar\alpha} \Biggl\{\sum_{k=1}^n
[h_{\alpha^\circ}(k)-h_\alpha(k)] b(k) \xi(k)
\\
&&\qquad\hspace*{17.77pt}{}
-z\Biggl[\sum_{k=1}^n
[h_{\alpha^\circ}(k)-h_\alpha(k)]^2 b^2(k) \Biggr]^{q/2} \Biggr\}_+^p
\le\frac{C(p,q)}{z^{p/(q-1)}}.
\end{eqnarray*}
\end{lemma}

The next fact is essential for bounding cross terms in the empirical risk.
\begin{lemma}\label{lemma.5} Let $\alpha^\circ$ be a given smoothing
parameter. Then for any $p\in[1,2)$,
there exists a constant $C(p)$ such that for any data-driven smoothing
parameter $\hat\alpha$,
%
%
\begin{eqnarray}\label{equ.28}\qquad
&&\mathbf{E} \Biggl| \sum_{k=1}^\infty
[h_{\hat\alpha}(k)-h_{\alpha^\circ}(k)]\lambda^{-1/2}(k) \theta(k)
\xi(k) \Biggr|^p\nonumber\\
&&\qquad \le C(p)
\Bigl\{\mathbf{E}\max_{k} \lambda^{-1}(k)
h_{\hat\alpha}^2(k) \Bigr\}^{p/2} \Biggl\{ \sum_{k=1}^\infty
[1-h_{\alpha^\circ}(k)]^2\theta^2(k) \Biggr\}^{p/2}
\\
&&\qquad\quad{} + C(p) \Bigl\{\max_{k} \lambda^{-1}(k)
h_{\alpha^\circ}^2(k) \Bigr\}^{p/2} \Biggl\{\mathbf{E} \sum_{k=1}^\infty
[1-h_{\hat\alpha}(k)]^2\theta^2(k) \Biggr\}^{p/2}.\nonumber
\end{eqnarray}
\end{lemma}
\begin{pf}
From Lemmas \ref{lemma.4} and \ref{lemma.3}, it follows that
%
%
\begin{eqnarray}\label{equ.29}
&&\mathbf{E} \Biggl| \sum_{k=1}^\infty
[h_{\hat\alpha}(k)-h_{\alpha^\circ}(k)]\lambda^{-1/2}(k) \theta(k)
\xi(k) \Biggr|^p\nonumber\\[-8pt]\\[-8pt]
&&\qquad \le C(p)
\Biggl\{\mathbf{E} \sum_{k=1}^\infty
[h_{\alpha^\circ}(k)-h_{\hat\alpha}(k)]^2\lambda^{-1}(k)\theta^2(k)
\Biggr\}^{p/2}.\nonumber
\end{eqnarray}

To bound from above the right-hand side of (\ref{equ.19}), we
use that $h_\alpha(\cdot) $ is a family of ordered smoothers. With
this in mind, let us assume for definiteness that $h_{\alpha_1}(k)\ge
h_{\alpha_2}(k)$ for $\alpha_1\le\alpha_2$. Then,
if $\hat\alpha\ge\alpha^\circ$, we get
\[
\frac{h_{\hat\alpha}(k)}{h_{\alpha^\circ}(k)}\le1,\qquad \frac{h_{\hat
\alpha}(k)}{h_{\alpha^\circ}(k)}\ge h_{\hat\alpha}(k)
\]
and thus we obtain
%
%
\begin{eqnarray}\label{equ.30}
&&\sum_{k=1}^\infty
[h_{\alpha\circ}(k)-h_{\hat\alpha}(k)]^2\lambda^{-1}(k)\theta^2(k)
\nonumber\\
&&\qquad =\sum_{k=1}^\infty
H^2_{\alpha^\circ}(k)
\biggl[1-\frac{h_{\hat\alpha}(k)}{h_{\alpha^\circ}(k)} \biggr]^2\lambda^{-2}(k)
\theta^2(k)\\
&&\qquad
\le\max_k \lambda^{-1}(k) h_{\alpha^\circ}^2(k) \sum
_{k=1}^\infty
[1-h_{\hat\alpha}(k)]^2\theta^2(k).\nonumber
\end{eqnarray}
Similarly, if $\hat\alpha< \alpha^\circ$, then
%
%
\begin{eqnarray}\label{equ.31}
&&\sum_{k=1}^\infty
[h_{\alpha\circ}(k)-h_{\hat\alpha}(k)]^2\lambda^{-1}(k)\theta^2(k)
\nonumber\\
&&\qquad
\le\max_k \lambda^{-1}(k) h_{\hat\alpha}^2(k) \sum_{k=1}^\infty
[1-h_{\alpha^\circ}(k)]^2\theta^2(k)\\
&&\qquad\quad{} +\max_k \lambda^{-1}(k) h_{\alpha^\circ}^2(k) \sum_{k=1}^\infty
[1-h_{\hat\alpha}(k)]^2\theta^2(k).\nonumber
\end{eqnarray}

Therefore, combining (\ref{equ.30}) and (\ref{equ.31}), we get (\ref{equ.28}).
\end{pf}

The following important ordered process is defined by
\[
\zeta_\alpha=\sum_{k=1}^n h_{\alpha}(k)b(k)[\xi^2(k)-1],
\]
where $\xi(k)$ are i.i.d. $\mathcal{N}(0,1)$.
Let
\[
\sigma_\alpha= \Biggl[2\sum_{k=1}^n h_{\alpha}^2(k)b^2(k) \Biggr]^{1/2}.
\]
It is easy to check that
$
\min\{\mathbf{E}\zeta_{\alpha_1}^2,\mathbf{E}\zeta_{\alpha_2}^2
\}\le\mathbf{E}\zeta_{\alpha_2}\zeta_{\alpha_1}
$
and thus
\[
|\sigma_{\alpha_1}^2-\sigma_{\alpha_1}^2|\ge\|h_{\alpha
_1}-h_{\alpha_2}\|^2_b,
\]
where
\[
\|h_{\alpha_1}-h_{\alpha_2}\|^2_b =\sum_{k=1}^n b^2(k)
[h_{\alpha_1}(k)-h_{\alpha_2}(k)]^2.
\]
Hence, in order to apply Lemma \ref{lemma.2}, it remains to check that
for some $\Lambda>0$
\[
\sup_{\alpha_1,\alpha_2}\mathbf{E}\exp\biggl[\Lambda\frac{\Delta
_\zeta(\alpha_1,\alpha_2)}{\|h_{\alpha_1}-h_{\alpha_2}\|_b} \biggr]<
\infty.
\]

We have
%
%
\begin{eqnarray}\label{equ.32}
&& \mathbf{E}\exp\biggl[\Lambda\frac{\Delta_\eta(\alpha_1,\alpha_2)}{\|
h_{\alpha_1}-h_{\alpha_2}\|_b} \biggr]\nonumber\\
&&\qquad =
\exp\Biggl\{-\frac{\Lambda}{\sqrt{2}\|h_{\alpha_1}-h_{\alpha_2}\|
_b}\sum_{k=1}^n
b(k)[h_{\alpha_1}(k)-h_{\alpha_2}
(k)] \\
&&\qquad\quad\hspace*{20.58pt}{} -\frac{1}{2}\sum_{k=1}^n
\log\biggl[1-\sqrt{2}\Lambda\frac{b(k)[h_{\alpha_1}(k)-h_{\alpha
_2}(k)]}{\|h_{\alpha_1}-h_{\alpha_2}\|_b} \biggr]
\Biggr\}.\nonumber
\end{eqnarray}
Since obviously
\[
\max_k\{|b(k)||h_{\alpha_1}(k)-h_{\alpha_2}(k)|\}\le
\|h_{\alpha_1}-h_{\alpha_2}\|_b,
\]
then using the Taylor expansion for $\log(1-\cdot)$ at the
right-hand side of (\ref{equ.32}), we get for $\Lambda\le1/2$
\[
\mathbf{E}\exp\biggl(\Lambda\frac{\Delta_\eta(\alpha_1,\alpha_2)}{\|
h_{\alpha_1}-h_{\alpha_2}\|_b} \biggr)\le
\exp(C\lambda^2),
\]
thus proving (\ref{equ.23}). Hence, with the help of Lemma \ref
{lemma.2} we
obtain the following fact.
\begin{lemma} \label{lemma.6}
For any $z>0$,
\[
\mathbf{E}\sup_{\alpha\in(0,\bar{\alpha}]} \Biggl[\sum_{k=1}^n
h_\alpha(k)b(k)[\xi^2(k)-1]-z \Biggl[2\sum_{k=1}^n
h_\alpha^2(k)b^2(k) \Biggr]^{q/2} \Biggr]_+^p \le
\frac{C(p,q)}{z^{p/(q-1)}}.
\]
\end{lemma}

\subsection[Proof of Theorem 1]{Proof of Theorem \protect\ref{theorem.1}}\label{sec24}
The next lemma describes some basic properties of the universal penalty
defined by (\ref{equ.12})--(\ref{equ.14}).
\begin{lemma}\label{lemma.7} For any $\alpha\in(0,\bar\alpha]$,
%
%
\begin{eqnarray}\label{equ.33}\qquad\quad\qquad
\log\frac{D(\alpha)}{D(\bar{\alpha})}&\le&\mu_\alpha\frac
{Q^{\circ}(\alpha)}{D(\alpha)},\\
\label{equ.34}
\mu_\alpha&\ge&\min\Biggl\{\frac{1}{2}\sqrt{\log\frac{D(\alpha
)}{D(\bar{\alpha})}}, \frac14 \Biggr\},
\\
\label{equ.35}
\frac{Q^{\circ}(\alpha)}{D(\bar{\alpha})}&\ge&\frac{D(\alpha
)}{D(\bar{\alpha})} \biggl(\log\frac{D(\alpha)}{D(\bar{\alpha})} \biggr)^{1/2},
\\
\label{equ.36}
\frac{D(\alpha)}{D(\bar{\alpha})}&\ge&\frac{\mu_\alpha Q^{\circ
}(\alpha)}{D(\bar{\alpha})}\Big/\log\frac{\mu_\alpha Q^{\circ
}(\alpha)}{D(\bar{\alpha})}\qquad \mbox{if } {D(\alpha)}\ge\exp
(2)D(\bar{\alpha}),\\
\label{equ.37}
\frac{D(\alpha_1)}{D(\alpha_2)}&\le&\frac{Q^{\circ}(\alpha
_1)}{Q^{\circ}(\alpha_2)},\qquad \alpha_1 \le\alpha_2.
\end{eqnarray}
\end{lemma}
\begin{pf}
It follows from (\ref{equ.12})--(\ref{equ.14}) that
\[
\mu_\alpha\frac{Q^{\circ}(\alpha)}{D(\alpha)}+\sum_{k=1}^{n} \biggl\{
\frac{1}{2}\log[1-2\mu_\alpha\rho_\alpha(k)]+\mu_\alpha\rho
_\alpha(k) \biggr\}=\log\frac{D(\alpha)}{D(\bar{\alpha})}
\]
and together with the inequality $\log(1-2x)/2+x\le0$, we
get (\ref{equ.33}).

To verify (\ref{equ.34}), note that
\[
F(x)\le\frac{2x^2}{1-2x}
\]
and therefore, for any $\mu\in[0, 1/4]$,
\[
\sum_{k=1}^n F[\mu\rho_\alpha(k)]\le4\mu^2.
\]
So, if $\mu_\alpha\le1/4$, then
\[
\mu_\alpha^2 \ge\frac14 \sum_{k=1}^n F[\mu_\alpha\rho_\alpha
(k)]=\frac14\log\frac{D(\alpha)}{D(\bar{\alpha})},
\]
thus proving (\ref{equ.34}).

Next, note that the following inequality holds:
%
%
\begin{equation}\label{equ.38}
F(x)\ge\frac{x^2}{1-2x}
\end{equation}
since
\[
f(x)=F(x)-\frac{x^2}{1-2x}=\frac{1}{2}\log(1-2x)+x+\frac{x^2}{1-2x}
\]
is a nonnegative function for $x\ge0$ because
\[
f'(x)=\frac{2x^2}{(1-2x)^2}\ge0
\]
and $f(0)=0$.

According to (\ref{equ.38}), $F(x)\ge x^2$ and therefore, by
(\ref{equ.13}),
\[
\mu_\alpha^2 \le\frac{D(\alpha)}{D(\bar{\alpha})}.
\]
Substituting this inequality into (\ref{equ.33}), we get (\ref{equ.35}).

We now turn to the proof of (\ref{equ.36}). Again, combining (\ref
{equ.38}) with (\ref{equ.12})--(\ref{equ.14}), we arrive at
%
%
\begin{eqnarray}\label{equ.39}\quad
Q^{\circ}(\alpha)&=&2D(\alpha)\mu_\alpha\sum_{k=1}^n\frac{\rho
^2_\alpha(k)}{1-2\mu_\alpha
\rho_\alpha(k)}\le\frac{2 D(\alpha)}{\mu_\alpha}\sum_{k=1}^n
F[\mu_\alpha\rho_\alpha(k)]\nonumber\\[-8pt]\\[-8pt]
&=& \frac{2D(\alpha)}{\mu_\alpha} \log
\frac{D(\alpha)}{D(\bar{\alpha})}\nonumber
\end{eqnarray}
and to get (\ref{equ.36}) it remains to invert this equation. We
proceed to show
that if $x\ge\exp(2)$, then the inequality
$
x\log(x)\ge y
$
implies
%
%
\begin{equation}\label{equ.40}
x\ge\frac{y}{\log(y)}.
\end{equation}
It is clear that $G(y)=y/\log(y)$ is an increasing function when $y\ge
\exp(1)$ and (\ref{equ.40}) holds since
\[
x \ge\frac{x\log(x)}{\log(x\log(x))}=\frac{x}{1+\log\log
(x)/\log(x)}.
\]
Inverting (\ref{equ.39}) with the help of (\ref{equ.40}), we finish
the proof of (\ref{equ.36}).

Finally, (\ref{equ.37}) follows from the fact that
\[
g(\alpha)=\frac{Q^{\circ}(\alpha)}{D(\alpha)}=2\mu_\alpha\sum
_{k=1}^n\frac{\rho^2_\alpha(k)}{1-2\mu_\alpha\rho_\alpha(k)}
\]
is a decreasing function in $\alpha>0$. To check this, let us note
that $g(\alpha)$ is a root of the equation
\[
\exp\Bigl\{\inf_{\mu\ge0} [ -G_\alpha
(\mu)- \mu g(\alpha)] \Bigr\}=\frac{D(\bar{\alpha})}{D(\alpha)},
\]
where
\[
G_\alpha(\mu) =\sum_{k=1}^n \biggl[\mu\rho_\alpha(k)+\frac12\log
[1-2\mu\rho_\alpha(k)] \biggr].
\]
However $\inf_{\mu\ge0} [ -G_\alpha
(\mu)- \mu x]$ is obviously a decreasing function in $x$ and therefore
if $ {D(\bar{\alpha})}/{D(\alpha)}$ is decreasing in $\alpha$, then
$g(\alpha)$ is decreasing in $\alpha$ too.
\end{pf}

We are now in a position to prove Theorem \ref{theorem.1}.
Let $\alpha_k, k=0,\ldots,$ be the decreasing sequence
defined as follows:
\[
\alpha_0=\bar{\alpha},\qquad Q^\circ(\alpha_k)=(1+\delta)^{k-1}Q^\circ
(\alpha_1),
\]
where $\delta<1/2$ is a small positive number which will be chosen
later on, and $\alpha_1$ is a root of equation
\[
D(\alpha_1)=D(\bar{\alpha})\exp(2).
\]
Denote for brevity $D_k=D(\alpha_k)$ and $Q_k=Q^{\circ}(\alpha_k)$.

We begin with the simple inequality
\begin{eqnarray*}
&&\mathbf{E}\sup_{\alpha\le\bar{\alpha}} [\eta_\alpha-(1+\gamma
)Q^{\circ}(\alpha) ]_+^{1+r}\\
&&\qquad
\le\sum_{k=1}^n \mathbf{E}\sup_{\alpha_{k}\le\alpha< \alpha
_{k-1} } [\zeta_\alpha-(1+\gamma){Q}_{k-1} ]_+^{1+r}\\
&&\qquad = \sum_{k=1}^n \mathbf{E} \Bigl[\zeta_{\alpha_{k}}-\varepsilon\gamma
{Q}_{k-1}+\sup_{\alpha_{k}\le\alpha< \alpha_{k-1}
}[\zeta_\alpha-\zeta_{\alpha_{k}}]-(1+\gamma-\varepsilon\gamma
){Q}_{k-1} \Bigr]_+^{1+r}.
\end{eqnarray*}
Using that $[x+y]_+^{1+r}\le2^r[x]_+^{1+r}+2^r[y]_+^{1+r}$, we can
continue the above equation as follows:
%
%
\begin{eqnarray}\label{equ.41}
&&
\mathbf{E}\sup_{\alpha\le
\bar{\alpha}} [\eta_\alpha-(1+\gamma)Q^{\circ}(\alpha) ]_+^{1+r}
\nonumber\\
&&\qquad
\le2^r\sum_{k=1}^n \mathbf{E}
[\zeta_{\alpha_{k}}-(1+\gamma-\varepsilon\gamma){Q}_{k-1}
]_+^{1+r}\\
&&\qquad\quad{} +2^r\sum_{k=1}^n \mathbf{E} \Bigl[\sup_{\alpha_{k}\le\alpha<
\alpha_{k-1}}[\zeta_\alpha-\zeta_{\alpha_{k}}]-\varepsilon
\gamma{Q}_{k-1} \Bigr]_+^{1+r}.\nonumber
\end{eqnarray}

We control the first term
\[
\Delta_1(\gamma,\varepsilon)\stackrel{\mathrm{def}}{=}\sum_{k=1}^n
\mathbf
{E} [\zeta_{\alpha_{k}}-(1+\gamma-\varepsilon\gamma){Q}_{k-1} ]_+^{1+r}
\]
at the right-hand side of (\ref{equ.41}) with the help of (\ref
{equ.11}). Thus, we obtain for any $\tilde\lambda_k>0$
%
%
\begin{eqnarray}\label{equ.42}\quad
\Delta_1(\gamma,\varepsilon) &\le& \mathbf{E}|\zeta_{\alpha_{1}}|^{1+r}+
\sum_{k=2}^n {D}_{k}^{1+r}\mathbf{E} \biggl[\frac{\zeta_{\alpha
_{k}}}{{D}_{k}}-(1+\gamma-\varepsilon\gamma)
\frac{{Q}_{k-1}}{{D}_{k}} \biggr]_+^{1+r}\nonumber\\
&\le& C{D}^{1+r}(\bar{\alpha})\nonumber\\
&&{} + \Gamma(1+r)\sum_{k=2}^n
\tilde{\lambda}^{-1-r}_{k} {D}^{r}_{k}\exp\biggl\{-\tilde{\lambda}_{k}
\frac{{Q}_{k}}{{D}_{k}} \biggl[
\frac{\gamma(1-\varepsilon){Q}_{k-1}}{{Q}_{k}}
\\
&&\qquad\quad\hspace*{141.6pt}{} - \biggl(1-\frac{{Q}_{k-1}}{{Q}_k} \biggr) \biggr] \biggr\}\nonumber\\
&&\qquad\quad\hspace*{-23pt}{}\times {D}_{k}\mathbf{E}\exp\biggl(\tilde{\lambda}_{k}
\frac{\zeta_{\alpha_{k}}}{{D}_{k}}-\tilde{\lambda}_{k}
\frac{{Q}_{k}}{{D}_{k}} \biggr).\nonumber
\end{eqnarray}
According to (\ref{equ.12}) and (\ref{equ.13}), we have with $\tilde
{\lambda}_k =\mu_{\alpha_k}$
\[
D_{k}\mathbf{E}\exp\biggl(\tilde{\lambda}_{k} \frac{\zeta_{\alpha
_{k}}}{{D}_{k}}-\tilde{\lambda}_{k} \frac{{Q}_{k}}{{D}_{k}} \biggr)=D_0,
\]
and substituting this into (\ref{equ.42}), we get
%
%
\begin{eqnarray}\label{equ.43.1}\quad
\Delta_1(\gamma,\varepsilon)
&\le& C{D}^{1+r}_0
\Biggl\{1+\sum_{k=2}^n \tilde{\lambda}^{-1-r}_{k}
\exp\biggl[-\tilde{\lambda}_{k}
\frac{(\gamma(1-\varepsilon)-\delta){Q}_{k}}{(1+\delta){D}_{k}}\nonumber\\[-8pt]\\[-8pt]
&&\hspace*{175.5pt}{} +r\log
\frac{{D}_{k}}{D_0} \biggr] \Biggr\}.\nonumber
\end{eqnarray}
Since by (\ref{equ.33})
%
%
\begin{equation}\label{equ.33.x}
\tilde{\lambda}_{k}\frac{{Q}_{k}}{D_k}\ge\log\frac{{D}_{k}}{D_0},
\end{equation}
and according to (\ref{equ.34}), $\tilde{\lambda}_{k}$ is bounded
from below by a constant, we obtain from (\ref{equ.43.1})
%
%
\begin{equation}\label{equ.43.2}\quad
\Delta_1(\gamma,\varepsilon) \le C{D}^{1+r}_0 \Biggl\{1 +
\sum_{k=2}^n \exp\biggl[- \biggl(\frac{\gamma(1-\varepsilon)-\delta}{1+\delta
}-r \biggr)\log
\frac{D_k}{D_0} \biggr] \Biggr\}.
\end{equation}
Next, according to (\ref{equ.36}) and (\ref{equ.34}), we get
%
%
\begin{equation}\label{equ.34.36}
\log\frac{D_k}{D_0}\ge\frac{CQ_k}{Q_0}\log^{-1} \frac{CQ_k}{Q_0}
\end{equation}
and with this inequality we obtain from (\ref{equ.43.2})
\begin{eqnarray*}
\Delta_1(\gamma,\varepsilon) &\le&
CD^{1+r}_0\sum_{k=2}^n \exp\biggl\{- \biggl(\frac{\gamma(1-\varepsilon)-\delta
}{1+\delta}-r \biggr)
\log\biggl[ \frac{CQ_k}{Q_0}\log^{-1} \frac{CQ_k}{Q_0} \biggr] \biggr\}\\
&\le&
CD^{1+r}_0\sum_{k=2}^n \exp\biggl\{- \biggl(\frac{\gamma(1-\varepsilon)-\delta
}{1+\delta}-r \biggr)
[k\log(1+\delta)-\log(k\delta)] \biggr\}.
\end{eqnarray*}
Finally, one can check by the Laplace method that
\[
\sum_{k=1}^\infty\exp[-z_1 k+z_2 \log(k)]\le\frac{C}{z_1} \biggl(\frac
{z_2}{z_1} \biggr)^{z_2},\qquad z_1,z_2> 0,
\]
thus yielding
%
%
\begin{equation}\label{equ.43}
\Delta_1(\gamma,\varepsilon)\le\frac{CD^{1+r}_0}{\delta[\gamma
(1-\varepsilon)-\delta-(1+\delta)r]_+ }.
\end{equation}

Our next step is to bound from above the last term in (\ref{equ.41}), namely,
\[
\Delta_2(\gamma,\varepsilon)\stackrel{\mathrm{def}}{=}\sum_{k=1}^n
\mathbf
{E} \Bigl[\sup_{\alpha_{k}\le\alpha< \alpha_{k-1}}[\zeta_\alpha-\zeta
_{\alpha_{k}}]-\varepsilon\gamma{Q}_{k-1} \Bigr]_+^{1+r}.
\]
Consider the following random processes:
\[
\tilde\zeta_\alpha(k)=\zeta_{\alpha}-\zeta_{\alpha_k},\qquad t\in
[0,\alpha_{k-1}-\alpha_{k}].
\]
Denote for brevity
$
\tilde{\sigma}_\alpha(k)=\sqrt{\mathbf{E}\tilde\zeta_\alpha^2(k)}
$.
Noticing that $\tilde h_{\alpha}(k) =2h_{\alpha}(k)-h_{\alpha}^2(k)$
is a family of ordered smoothers, it is easy to check that
\[
\tilde\sigma_u^2(k)-\tilde\sigma_v^2(k) \ge\mathbf{E}[\tilde
{\zeta}_u(k)-\tilde{\zeta}_v(k)]^2,\qquad u\le v.
\]
According to the Taylor formula, for all
$u\ge v$ and all $\lambda\ge0$,
\[
\mathbf{E}\exp\biggl\{\lambda\frac{\tilde{\zeta}_u(k)-\tilde{\zeta
}_v(k)}{\mathbf{E}^{1/2}[\tilde{\zeta}_u(k)-\tilde{\zeta}_v(k)]^2}
\biggr\}\le\exp\biggl(\frac{\lambda^2}{2} \biggr),
\]
and applying Lemma \ref{lemma.1} and (\ref{equ.11}), we obtain for
any $\lambda_k\ge0$
\begin{eqnarray*}
\Delta_2(\gamma,\varepsilon)&=&\sum_{k=1}^n \mathbf{E} \Bigl[\sup_{\alpha
\in[0,\alpha_{k-1}-
\alpha_{k}]}\tilde{\zeta}_\alpha(k)-\varepsilon\gamma{Q}_{k-1} \Bigr]_+^{1+r}
\\
&\le& CD^{1+r}_0
+ C\sum_{k=2}^n\tilde\sigma_{\alpha_{k-1}-\alpha
_{k}}^{1+r}(k)\lambda_k^{-1-r}\\
&&\hspace*{69pt}{}\times
\exp\biggl[- \frac{ \lambda_k\varepsilon\gamma{Q}_{k-1}}{\tilde\sigma
_{\alpha_{k-1}-\alpha_{k}}(k)}
+\frac{\lambda_k^2}{(\sqrt{2}-1)^2}
\biggr].
\end{eqnarray*}
Substituting
\[
\lambda_k=\frac{(\sqrt{2}-1)^2 \varepsilon\gamma{Q}_{k-1}}{2\tilde
\sigma_{\alpha_{k-1}-\alpha_{k}}(k)}
\]
into the above equation
and noticing that $ {Q}_k \ge{D}_k$, we obtain
%
%
\begin{eqnarray}\label{equ.44.1}\quad
\Delta_2(\gamma,\varepsilon) &\le& CD^{1+r}_0
+ \frac{C}{(\varepsilon\gamma)^{1+r}}\sum_{k=2}^n\tilde\sigma
_{\alpha_{k-1}-\alpha_{k}}^{2(1+r)}(k)
{Q}^{-1-r}_{k-1}\nonumber\\
&&\hspace*{99pt}{} \times \exp\biggl[- \frac{ (\sqrt{2}-1)^2\varepsilon^2\gamma
^2{Q}^2_{k-1}}{4\tilde\sigma_{\alpha_{k-1}-\alpha_{k}}^2(k)} \biggr]
\nonumber\\[-8pt]\\[-8pt]
&\le& CD^{1+r}_0+ \frac{C}{(\varepsilon\gamma)^{1+r}}\sum
_{k=2}^n[{D}^2_k-{D}^2_{k-1}]^{1+r}
{Q}^{-1-r}_{k-1} \nonumber\\
&&\hspace*{98.3pt}{}\times \exp\biggl[- \frac{ (\sqrt{2}-1)^2\varepsilon^2\gamma
^2{Q}^2_{k-1}}{4({D}^2_k-{D}^2_{k-1})} \biggr].\nonumber
\end{eqnarray}
According to (\ref{equ.37}),
\[
\frac{D_k^2}{D_{k-1}^2}\le\frac{Q_{k}^2}{Q_{k-1}^2}\le(1+\delta)^2
\]
and with this inequality we continue (\ref{equ.44.1}) as follows:
\begin{eqnarray*}
\Delta_2(\gamma,\varepsilon)&\le& CD^{1+r}_0 + \frac{C\delta
^{1+r}}{(\varepsilon\gamma)^{1+r}}\sum_{k=2}^n \biggl(\frac{{D}_{k-1}}{{Q}_{k-1}}
\biggr)^{1+r}{D}_{k-1}^{1+r} \\
&&\hspace*{99pt}{}\times\exp\biggl[- \frac{ (\sqrt{2}-1)^2\varepsilon^2\gamma
^2{Q}_{k-1}^2}{8\delta{D}_{k-1}^2
} \biggr].
\end{eqnarray*}
Next, substituting (\ref{equ.33.x}) and (\ref{equ.34.36}) into this
equation, we get
\begin{eqnarray*}
\Delta_2(\gamma,\varepsilon)&\le& CD^{1+r}_0
 + \frac{C\delta^{1+r}}{(\varepsilon\gamma)^{1+r}}\sum_{k=2}^n
{Q}_{k-1}^{1+r} \exp\biggl[- \frac{ (\sqrt{2}-1)^2\varepsilon^2\gamma
^2}{8\delta\tilde{\mu}_{\alpha_{k-1}}^2
}\log^2 \frac{{D}_{k-1}}{{D}_0} \biggr]
\\
&\le& CD^{1+r}_0+ \frac{C\delta^{1+r}}{(\varepsilon\gamma
)^{1+r}}\sum_{k=2}^n
{Q}_{k-1}^{1+r} \exp\biggl[- \frac{ C\varepsilon^2\gamma^2}{\delta}\log^2
\frac{{Q}_{k-1}}{{D}_0} \biggr]
\\
&\le&
\frac{C\delta^{1+r}D^{1+r}_0}{(\varepsilon\gamma)^{1+r}}\sum_{k=2}^n
\exp\biggl[\delta(r+1) k-\frac{C\varepsilon^2\gamma^2 }{\delta}[k\log
(1+\delta)-\log(\delta k)]^2 \biggr].
\end{eqnarray*}
Bounding the last sum in this display with the help of the Laplace method,
we get
\begin{eqnarray*}
&&\sum_{k=2}^n
\exp\biggl\{\delta(r+1) k-\frac{C\varepsilon^2\gamma^2 }{\delta}[k\log
(1+\delta)-\log(\delta k)]^2 \biggr\}\\
&&\qquad\le\exp\biggl[\frac{C\delta
^3}{\varepsilon^2\gamma^2\log^2(1+\delta)}
\biggr]\frac{C\sqrt{\delta}}{\varepsilon\gamma\log(1+\delta)},
\end{eqnarray*}
and therefore with $\delta=\varepsilon^2\gamma^2$ we obtain
%
%
\begin{equation}\label{equ.44}
\Delta_2(\gamma,\varepsilon) \le\frac{CD^{1+r}_0}{(\varepsilon\gamma)^{1-r}}.
\end{equation}
With this $\delta$, equation (\ref{equ.43}) becomes
%
%
\begin{equation}\label{equ.45}
\Delta_1(\gamma,\varepsilon)\le\frac{CD^{1+r}(\bar{\alpha
})}{(\gamma\varepsilon)^2[\gamma-r-\gamma\varepsilon-(1+r)(\gamma
\varepsilon)^2]_+}.
\end{equation}

Let $\varepsilon$ be a positive root of the equation
\[
\gamma-r-\gamma\varepsilon-(1+r)(\gamma\varepsilon)^2=\gamma\varepsilon,
\]
that is,
\[
\varepsilon=-\frac{1}{(r+1)\gamma}+\frac{1}{\gamma}\sqrt{\frac
{1}{(1+r)^2}+\frac{\gamma-r}{1+r}}.
\]
Substituting this $\varepsilon$ into (\ref{equ.44}) and (\ref{equ.45}) and
combining the obtained inequalities with (\ref{equ.41}), we finish
the proof.

\subsection[Proof of Theorem 2]{Proof of Theorem \protect\ref{theorem.2}}\label{sec25}

The first step in the proof of this theorem is to show that the
data-driven parameter $\hat\alpha$ defined by (\ref{equ.4}) cannot
be very small, or equivalently, that the ratio
$D(\hat\alpha)/D(\bar{\alpha})$ is not large.
\begin{lemma}\label{lemma.8} For any given $\alpha^\circ\le\bar
{\alpha}$ and $\gamma>0$ the following upper bound holds
%
%
\begin{equation}\label{equ.46}\quad
\biggl\{\mathbf{E} \biggl[\frac{D(\hat\alpha)}{D(\bar{\alpha})} \biggr]^{1+\gamma
/4} \biggr\}^{1/(1+\gamma/4)}\le
\frac{C\bar{R}_{\alpha^\circ}(\theta)}{\sigma^2\gamma D(\bar
{\alpha})}
\log^{-1/2}\frac{C\bar{R}_{\alpha^\circ}(\theta)}{\sigma^2D(\bar
{\alpha})}
+ \frac{C}{\gamma^4}.
\end{equation}
\end{lemma}
\begin{pf}
According to the definition of the empirical risk minimization,
for any given $\alpha^\circ$, $R_{\hat\alpha}[Y,\operatorname{Pen}]\le R_{\alpha
^\circ}[Y,\operatorname{Pen}]$. One can check with a little algebra that this
inequality is equivalent
to [see (\ref{equ.7})]
%
%
\begin{eqnarray}\label{equ.47}\quad
&&\sum_{k=1}^n[1-h_{\hat\alpha}(k)]^2\langle
\theta,\psi_k\rangle^2 +\sigma^2\sum_{k=1}^n \lambda^{-1}(k)
h_{\hat\alpha}^2(k)\nonumber\\
&&\quad{} -\sigma^2 \sum_{k=1}^n \lambda^{-1}(k)\tilde{h}_{\hat\alpha
}(k)[\xi^2(k)-1]+(1+\gamma)
\sigma^2Q^{\circ}(\hat\alpha)\nonumber\\
&&\quad{} +2\sigma\sum_{k=1}^n \lambda^{-1/2}(k)[1-h_{\hat\alpha
}(k)]^2\xi(k)\theta(k)
\nonumber\\[-8pt]\\[-8pt]
&&\qquad\le\sum_{k=1}^n[1-h_{\alpha^\circ}(k)]^2\langle
\theta,\psi_k\rangle^2 +\sigma^2\sum_{k=1}^n \lambda^{-1}(k)
h_{\alpha^\circ}^2(k)\nonumber\\
&&\qquad\quad{} -\sigma^2 \sum_{k=1}^n \lambda^{-1}(k)\tilde{h}_{\alpha^\circ
}(k)[\xi^2(k)-1]+(1+\gamma)
\sigma^2Q^{\circ}(\alpha^\circ)\nonumber\\
&&\qquad\quad{} +2\sigma\sum_{k=1}^n \lambda^{-1/2}(k)[1-h_{\alpha^\circ
}(k)]^2\xi(k)\theta(k),\nonumber
\end{eqnarray}
where
$
\tilde{h}_\alpha(k)=2h_\alpha(k)-h^2_\alpha(k)
$. Next, representing
\[
(1+\gamma)Q^{\circ}(\hat\alpha)= \biggl(1+\frac{\gamma}{2} \biggr)Q^{\circ
}(\hat\alpha)
+\frac{\gamma}{2} Q^{\circ}(\hat\alpha),
\]
we obtain from (\ref{equ.47})
%
%
\begin{eqnarray}\label{equ.48}\qquad
\frac{\gamma\sigma^2}{2}Q^{\circ}(\hat\alpha)
&\le&
\bar{R}_{\alpha^\circ}(\theta) +\sigma^2 \sum_{k=1}^n \lambda
^{-1}(k)\tilde{h}_{\alpha^\circ}(k)[\xi^2(k)-1]\nonumber\\
&&{} + \sigma^2\sup_{\alpha\le\bar{\alpha}} \Biggl[\sum_{k=1}^n \lambda
^{-1}(k)\tilde{h}_{\alpha}(k)[\xi^2(k)-1]
- \biggl(1+\frac{\gamma}{2} \biggr)
Q^{\circ}(\alpha) \Biggr]_+\\
&&{} + 2\sigma\sum_{k=1}^n \lambda^{-1/2}(k)[\tilde{h}_{\hat\alpha
}(k)-\tilde{h}_{\alpha^\circ}(k)
]\xi(k)\theta(k) - L_{\hat\alpha}(\theta).\nonumber
\end{eqnarray}

Since $\alpha^\circ$ is fixed, we get by Jensen's inequality
%
%
\begin{eqnarray}\label{equ.49}\qquad\quad
\mathbf{E} \Biggl|\sum_{k=1}^n \lambda^{-1}(k)\tilde{h}_{\alpha^\circ
}(k)[\xi^2(k)-1] \Biggr|^{1+\gamma/4}
&\le&
C \Biggl[\sum_{k=1}^n \lambda^{-2}(k)\tilde{h}_{\alpha^\circ}^2(k)
\Biggr]^{1/2+\gamma/8}\nonumber\\
&=& C[D(\alpha^\circ)]^{1+\gamma/4}\\
&\le& C[\sigma^{-2}\bar{R}_{\alpha^\circ}(\theta)]^{1+\gamma/4}.\nonumber
\end{eqnarray}

Next, by Theorem \ref{theorem.1},
%
%
\begin{eqnarray}\label{equ.50}
&&
\mathbf{E}\sup_{\alpha\le\bar{\alpha}} \Biggl[\sum_{k=1}^n \lambda
^{-1}(k)\tilde{h}_{\alpha}(k)[\xi^2(k)-1]
- \biggl(1+\frac{\gamma}{2} \biggr)
Q^{\circ}(\alpha) \Biggr]_+^{1+\gamma/4}\nonumber\\[-8pt]\\[-8pt]
&&\qquad
\le\frac{CD^{1+\gamma/4}(\bar{\alpha})}{\gamma^3}.\nonumber
\end{eqnarray}

The upper bound for the last line in (\ref{equ.48}) is a little bit
more tricky. Noticing that $\tilde{h}_{\alpha}(\cdot)$ is a family
of ordered smoothers, we get by Lemma \ref{lemma.4} that, for any
$\varepsilon>0$ and given $p\in(1,2)$,
%
%
\begin{eqnarray}\label{equ.51}
&&\mathbf{E} \Biggl|2\sigma\sum_{k=1}^n \lambda^{-1/2}(k)[\tilde{h}_{\hat
\alpha}(k)-\tilde{h}_{\alpha^\circ}(k)
]\xi(k)\theta(k)\nonumber\\
&&\quad{} - \varepsilon\Biggl[4\sigma^2\sum_{k=1}^n \lambda^{-1}(k)[\tilde{h}_{\hat
\alpha}(k)-\tilde{h}_{\alpha^\circ}(k)
]^2\theta^2(k) \Biggr]^{p/2} \Biggr|^{1+\gamma/4}\\
&&\qquad \le\frac{C(p)}{\varepsilon^{(1+\gamma/4)/(p-1)}}.\nonumber
\end{eqnarray}

To continue this inequality, note that
if $\hat\alpha\ge\alpha^\circ$, then
\[
\frac{\tilde{h}_{\hat\alpha}(k)}{\tilde{h}_{\alpha^\circ}(k)}\le
1,\qquad \frac{\tilde{h}_{\hat\alpha}(k)}{\tilde{h}_{\alpha^\circ
}(k)}\ge\tilde{h}_{\hat\alpha}(k)
\]
and therefore
%
%
\begin{eqnarray}\label{equ.52}
&&\sum_{k=1}^\infty
[\tilde{h}_{\alpha\circ}(k)-\tilde{h}_{\hat\alpha}(k)]^2\lambda
^{-1}(k)\theta^2(k)
\nonumber\\
&&\qquad =\sum_{k=1}^\infty
\tilde{h}^2_{\alpha^\circ}(k)
\biggl[1-\frac{\tilde{h}_{\hat\alpha}(k)}{\tilde{h}_{\alpha^\circ}(k)}
\biggr]^2\lambda^{-2}(k)
\theta^2(k)\nonumber\\[-8pt]\\[-8pt]
&&\qquad
\le\max_k \lambda^{-1}(k) \tilde{h}_{\alpha^\circ}^2(k) \sum
_{k=1}^\infty
[1-\tilde{h}_{\hat\alpha}(k)]^2\theta^2(k)
\nonumber\\
&&\qquad
\le4 \max_k \lambda^{-1}(k) {h}_{\alpha^\circ}^2(k) \sum
_{k=1}^\infty
[1-h_{\hat\alpha}(k)]^2\theta^2(k).\nonumber
\end{eqnarray}
Analogously, if $\hat\alpha< \alpha^\circ$, then
%
%
\begin{eqnarray}\label{equ.53}
&&\sum_{k=1}^n
[\tilde{h}_{\alpha\circ}(k)-\tilde{h}_{\hat\alpha}(k)]^2\lambda
^{-1}(k)\theta^2(k)
\nonumber\\
&&\qquad
\le\max_k \lambda^{-1}(k) \tilde{h}_{\hat\alpha}^2(k) \sum_{k=1}^n
[1-\tilde{h}_{\alpha^\circ}(k)]^2\theta^2(k)\\
&&\qquad \le4\max_k \lambda^{-1}(k) {h}_{\hat\alpha}^2(k) \sum_{k=1}^n
[1-{h}_{\alpha^\circ}(k)]^2\theta^2(k).\nonumber
\end{eqnarray}

Next, combining (\ref{equ.51})--(\ref{equ.53}) with Young's inequality,
%
%
\begin{equation}\label{equ.54}\quad
yx^q-x \le y^{-1/(q-1)} \bigl[q^{-q/(q-1)}-q^{1/(q-1)} \bigr],\qquad
x,y\ge0, q<1,
\end{equation}
gives
\begin{eqnarray*}
&&\mathbf{E} \Biggl|2\sigma\sum_{k=1}^n \lambda^{-1/2}(k)[\tilde
{h}_{\hat\alpha}(k)-\tilde{h}_{\alpha^\circ}(k)
]\xi(k)\theta(k) - L_{\hat\alpha}(\theta) \Biggr|^{1+\gamma/4}\\
&&\qquad\le C\mathbf{E} \Biggl|2\sigma\sum_{k=1}^n \lambda^{-1/2}(k)[\tilde
{h}_{\hat\alpha}(k)-\tilde{h}_{\alpha^\circ}(k)
]\xi(k)\theta(k)\\
&&\qquad\quad\hspace*{15.9pt}{} - \varepsilon\Biggl[4\sigma^2\sum_{k=1}^n \lambda^{-1 }(k)[\tilde
{h}_{\hat\alpha}(k)-\tilde{h}_{\alpha^\circ}(k)
]^2\theta^2(k) \Biggr]^{p/2} \Biggr|^{1+\gamma/4}\\
&&\qquad\quad{}
+ C \mathbf{E} \Biggl| \varepsilon\Biggl[4\sigma^2\sum_{k=1}^n {\lambda
^{-1}(k)}[\tilde{h}_{\hat\alpha}(k)-\tilde{h}_{\alpha^\circ}(k)
]^2\theta^2(k) \Biggr]^{p/2} - L_{\hat\alpha}(\theta) \Biggr|^{1+\gamma/4}
\\
&&\qquad \le\frac{C}{\varepsilon^{(1+\gamma/4)/(p-1)}} +\frac{C}{\varepsilon
^{2(1+\gamma/4)/(p-2)}} \Bigl[\sigma^2\max_k \lambda^{-1}(k) {h}_{\alpha
^\circ}^2(k) \Bigr]^{{p(1+\gamma/4)}/({2-p})}
\\
&&\qquad\quad{} +\frac{C}{\varepsilon^{2(1+\gamma/4)/(p-2)}} \Biggl\{\sum_{k=1}^\infty
[1-{h}_{\alpha^\circ}(k)]^2\theta^2(k) \Biggr\}^{{p(1+\gamma/4)}/({2-p})}.
\end{eqnarray*}
Therefore, minimizing the right-hand side at the above equation in
$\varepsilon>0$, we get
\begin{eqnarray*}
&&\mathbf{E} \Biggl|2\sigma\sum_{k=1}^n \lambda^{-1}(k)[\tilde{h}_{\hat
\alpha}(k)-\tilde{h}_{\alpha^\circ}(k)
]\xi(k)\theta(k) - L_{\hat\alpha}(\theta) \Biggr|^{1+\gamma/4}\\
&&\qquad \le C \Biggl\{\sum_{k=1}^\infty
[1-{h}_{\alpha^\circ}(k)]^2\theta^2(k)+\sigma^2\max_k \lambda
^{-1}(k) {h}_{\alpha^\circ}^2(k) \Biggr\}^{1+\gamma/4}.
\end{eqnarray*}
This equation and (\ref{equ.48})--(\ref{equ.50}) imply
\[
\gamma^{1+\gamma/4}\mathbf{E}[\sigma^2Q^{\circ}(\hat\alpha
)]^{1+\gamma/4}
\le C\bar{R}_{\alpha^\circ}^{1+\gamma/4}(\theta)+ \frac{C[\sigma
^2D(\bar{\alpha})]^{1+\gamma/4}}{\gamma^3}
\]
and by (\ref{equ.35}) we get
%
%
\begin{equation}\label{equ.55}\quad
\gamma^{1+\gamma/4}\mathbf{E} \biggl[\frac{D(\hat\alpha)}{D(\bar
{\alpha})}\log^{1/2}\frac{D(\hat\alpha)}{D(\bar{\alpha})}
\biggr]^{1+\gamma/4}
\le C \biggl[\frac{\bar{R}_{\alpha^\circ}(\theta)}{\sigma^2 D(\bar
{\alpha})} \biggr]^{1+\gamma/4}+ \frac{C}{\gamma^3}.
\end{equation}
It is easily seen that
%
%
\begin{eqnarray}\label{equ.56}\qquad
&&\mathbf{E} \biggl[\frac{D(\hat\alpha)}{D(\bar{\alpha})}\log
^{1/2}\frac
{D(\hat\alpha)}{D(\bar{\alpha})} \biggr]^{1+\gamma/4}\nonumber\\[-8pt]\\[-8pt]
&&\qquad =\frac{1}{(1+\gamma/4)^{1/2+\gamma/8}} \mathbf{E} \biggl[\frac
{D(\hat\alpha)}{D(\bar{\alpha})} \biggr]^{1+\gamma/4} \biggl[\log\biggl(\frac
{D(\hat\alpha)}{D(\bar{\alpha})} \biggr)^{1+\gamma/4}
\biggr]^{1/2+\gamma/8}.\nonumber
\end{eqnarray}

To finish the proof, let us consider the function $f(x)=x\log
^{1/2+\gamma/8}(x)$, $ x\ge1$. Computing its second order derivative,
one can easily check that $f(x)$ is convex for all $x\ge\exp(1)=\mathrm{e}$.
So, $f(x+\mathrm{e}-1)$ is convex for $x\ge1$. Note also that
there exists a constant $C>0$ such that for all $x\ge1$,
\[
f(x)\ge\tfrac{1}{2} f(x+\mathrm{e}-1)-C.
\]
Therefore according to (\ref{equ.56}) and Jensen's inequality,
\begin{eqnarray*}
&&\mathbf{E} \biggl(\frac{D(\hat\alpha)}{D(\bar{\alpha})}\log
^{1/2}\frac{D(\hat\alpha)}{D(\bar{\alpha})} \biggr)^{1+\gamma/4}\\
&&\qquad\ge C
\biggl[\mathbf{E} \biggl(\frac{D(\hat\alpha)}{D(\bar{\alpha})} \biggr)^{1+\gamma
/4}+\mathrm{e}-1 \biggr]
\\
&&\qquad\quad{} \times\biggl\{\log\biggl[\mathbf{E} \biggl(\frac{D(\hat\alpha)}{D(\bar{\alpha
})} \biggr)^{1+\gamma/4}+\mathrm{e}-1 \biggr] \biggr\}^{1/2+\gamma/8}-C.
\end{eqnarray*}

Finally, substituting this inequality into (\ref{equ.55}) and
inverting $f(x)$, we arrive at (\ref{equ.46}).
\end{pf}

The next lemma controls the cross term in the empirical risk.
\begin{lemma}\label{lemma.9} Let $\tilde{h}_\alpha^\varepsilon
(k)=[(1+2\varepsilon) h_{\alpha}(k)-\varepsilon h_{\alpha
}^2(k)]/(1+\varepsilon)$. Then for any given $\varepsilon\ge0$ and
$\alpha
^\circ\in(0,\bar\alpha]$,
%
%
\begin{eqnarray}\label{equ.59}\quad
&&2\sigma\mathbf{E} \Biggl|
\sum_{k=1}^n[1-\tilde{h}_{\hat\alpha}^\varepsilon(k)
]\theta(k)\lambda^{-1/2}(k)
\xi(k) \Biggr|\nonumber\\
&&\qquad
\le\biggl[\frac{C\bar{R}
_{\alpha{^\circ}}(\theta)}{\gamma}\log^{-1/2}\frac{C\bar{R}
_{\alpha{^\circ}}(\theta)}{\sigma^2D(\bar\alpha)}+\frac{C\sigma
^2D(\bar\alpha)}{\gamma^4} \biggr]^{1/2}\\
&&\qquad\quad{} \times\Biggl[
\mathbf{E}\sum_{k=1}^n[1-h_{\hat\alpha}(k)
]^2\theta^2(k)
+\sum_{k=1}^n[1-h_{\alpha^\circ}(k)
]^2\theta^2(k) \Biggr]^{1/2}.\nonumber
\end{eqnarray}
\end{lemma}
\begin{pf}
Since $\tilde{h}_\alpha^\varepsilon(k)$
is a family of ordered smoothers, combining Lemma \ref{lemma.5} with
the obvious inequalities $\max_k \lambda^{-1}(k) {h}_{ \alpha
}^2(k)\le D(\alpha)$ and $\tilde{h}_\alpha^\varepsilon(k)\ge h_\alpha
(k)$, we obtain
%
%
\begin{eqnarray}\label{equ.62}
&&
2\sigma\mathbf{E} \Biggl|
\sum_{k=1}^n[1-\tilde{h}_{\hat\alpha}^\varepsilon(k)
]\theta(k)\lambda^{-1/2}(k)
\xi(k) \Biggr|\nonumber\\
&&\qquad =2\sigma\mathbf{E} \Biggl|
\sum_{k=1}^n[h_{\alpha^\circ}^\varepsilon(k)-h_{\hat\alpha}^\varepsilon
(k)]\theta(k)
\lambda^{-1/2}(k)
\xi(k) \Biggr|
\nonumber\\[-8pt]\\[-8pt]
&&\qquad \le C \sigma\Biggl[\mathbf{E}D(\hat\alpha) \sum_{k=1}^n
[1-{h}_{\alpha^\circ}(k)]^2\theta^2(k) \Biggr]^{1/2}\nonumber\\
&&\qquad\quad{} +C\sigma\Biggl[D( \alpha^\circ) \mathbf{E} \sum_{k=1}^n
[1-{h}_{\hat\alpha}(k)]^2\theta^2(k) \Biggr]^{1/2}.\nonumber
\end{eqnarray}

Next, according to (\ref{equ.35}),
$
Q^{\circ}(\alpha)\ge D(\alpha)\sqrt{\log[D(\alpha)/D(\bar{\alpha})]}
$, and
we get
\[
D(\alpha^\circ)\le C\sigma^{-2} \bar{R}_{\alpha^\circ}(\theta
)\log^{-1/2}
\frac{\bar{R}_{\alpha^\circ}(\theta)}{\sigma^2D(\bar\alpha)}.
\]
Substituting this inequality and (\ref{equ.46}) in (\ref{equ.62}), we
obtain (\ref{equ.59}).
\end{pf}

We are now in a position to prove Theorem \ref{theorem.2}. Let
$\varepsilon\in(0,1]$ be a given number to be defined later on.
According to (\ref{equ.7}) and (\ref{equ.8}), we obtain the following
equation for the skewed excess risk:
%
%
\begin{eqnarray}\label{equ.57}\quad
\mathcal{E}(\varepsilon)&\stackrel{\mathrm{def}}{=}& \sup_{\theta\in
\mathbb{R}^n}\mathbf{E} \bigl\{\|\theta-\hat
\theta_{\hat\alpha}\|^2 -(1+\varepsilon) \{{R}_{\hat\alpha
}[Y,\operatorname{Pen}]+\mathcal{C}\} \bigr\}\nonumber\\
& = &\sup_{\theta\in
\mathbb{R}^n}\mathbf{E} \Biggl\{ -\varepsilon
\sum_{k=1}^n[1-h_{\hat\alpha}(k)]^2 \theta^2(k)-\varepsilon\sigma^2
\sum_{k=1}^n
\lambda^{-1}(k) h_{\hat\alpha}^2(k)\nonumber\\
&&\hspace*{33.4pt}{} - (1+\varepsilon)(1+\gamma)\sigma^2 Q^{\circ}(\hat\alpha)
\nonumber\\[-8pt]\\[-8pt]
&&\hspace*{33.4pt}{} -2\sigma
\sum_{k=1}^n\{1+\varepsilon-[(1+2\varepsilon)h_{\hat\alpha
}(k)-\varepsilon
h_{\hat\alpha}^2(k)
]\}\nonumber\\
&&\hspace*{72.8pt}{}\times\theta(k)\lambda^{-1/2}(k)
\xi(k) \nonumber\\
&&\hspace*{33.4pt}{} +\sigma^2\sum_{k=1}^n
\lambda^{-1}(k)[2(1+\varepsilon) h_{\hat\alpha}(k)-\varepsilon h_{\hat
\alpha}^2(k)][\xi^2(k)-1] \Biggr\}.\nonumber
\end{eqnarray}
To control the last line at the right-hand side of this equation, we
use that $h_\alpha^\varepsilon(k)=[2(1+\varepsilon) h_{\alpha
}(k)-\varepsilon h_{\alpha}^2(k)]/(2+\varepsilon)$
is a family of ordered smoothers. Hence, Lemmas \ref{lemma.3}, \ref
{lemma.6} and \ref{lemma.8} imply
%
%
\begin{eqnarray}\label{equ.58}
&&\sigma^2\mathbf{E}\sum_{k=1}^n
\lambda^{-1}(k) [2(1+\varepsilon) h_{\hat\alpha}(k)-\varepsilon h_{\hat
\alpha}^2(k)][\xi^2(k)-1]\nonumber\\[-8pt]\\[-8pt]
&&\qquad\le
C \frac{\bar{R}_{\alpha^\circ}(\theta)}{\gamma}
\log^{-1/2}\frac{C\bar{R}_{\alpha^\circ}(\theta)}{\sigma^2\gamma
D(\bar{\alpha})}+\frac{C\sigma^2D(\bar{\alpha})}{\gamma^4}.\nonumber
\end{eqnarray}

Next, substituting (\ref{equ.58}) and (\ref{equ.59}) into (\ref
{equ.57}), we obtain the following upper bound for the skewed excess
risk:
\[
\mathcal{E}(\varepsilon)\le\varepsilon\bar{R}_{\alpha^\circ}(\theta) +
\frac{C}{\varepsilon} \biggl[\frac{\bar{R}_{\alpha^\circ}(\theta)}{\gamma}
\log^{-1/2}\frac{C\bar{R}_{\alpha^\circ}(\theta)}{\sigma^2D(\bar
{\alpha})}+\frac{
\sigma^2D(\bar{\alpha})}{\gamma^4} \biggr].
\]

Finally, substituting this upper bound into
\[
\mathbf{E}\|\theta-\hat\theta_{\hat\alpha}\|^2 \le(1+\varepsilon
)\bar{R}_{\alpha^\circ}(\theta)+\mathcal{E}(\varepsilon)
\]
and minimizing the obtained inequality in $\varepsilon$, we get
\begin{eqnarray*}
\mathbf{E}\|\theta-\hat\theta_{\hat\alpha}\|^2 &\le& r(\theta) +
Cr(\theta)\inf_{\varepsilon\ge0} \biggl\{ \varepsilon+\frac{1 }{
\varepsilon} \biggl[ \frac{1}{\gamma}\log^{-1/2}\frac{Cr(\theta)}{\sigma^2
D(\bar{\alpha})}+\frac{C
\sigma^2D(\bar{\alpha})}{\gamma^4r(\theta)} \biggr] \biggr\}\\
&\le&
r(\theta) \biggl\{1 + \biggl[\frac{C}{\gamma}
\log^{-1/2}\frac{Cr(\theta)}{\sigma^2 D(\bar{\alpha})}+\frac{
C\sigma^2D(\bar{\alpha})}{\gamma^4r(\theta)} \biggr]^{1/2} \biggr\},
\end{eqnarray*}
thus finishing the proof.


\section*{Acknowledgments} The author wishes to thank two anonymous
referees for stimulating comments and remarks.

\printaddresses

\end{document}